\RequirePackage{fix-cm}
\documentclass[smallextended]{svjour3}       
\smartqed  
\usepackage{amsmath,amsfonts,latexsym,amssymb,amsbsy,tabularx}
\usepackage{mathtools,bm}
\usepackage{mathptmx}      

\usepackage{tikz}
\usepackage{graphicx}
\usepackage{color, xcolor}
\usepackage{adjustbox}

\usetikzlibrary{arrows, automata, shapes, decorations.pathreplacing, calligraphy, decorations.markings, positioning}

\usepackage{multirow,tabularx}
\usepackage{booktabs}

\usepackage{caption} 
\usepackage{subcaption}

\usepackage{algorithm}
\usepackage[noend]{algpseudocode}

\usepackage[inline]{enumitem}

\usepackage{url}
\usepackage[bookmarksnumbered=true]{hyperref}
\hypersetup{
	colorlinks,
	linkcolor={red!50!black},
	citecolor={blue!50!black},
	urlcolor={blue!80!black}
}

\usepackage{listings}

\definecolor{codegreen}{rgb}{0,0.6,0}
\definecolor{codegray}{rgb}{0.5,0.5,0.5}
\definecolor{codepurple}{rgb}{0.58,0,0.82}
\definecolor{backcolour}{rgb}{0.95,0.95,0.95}

\lstdefinestyle{mystyle}{
    backgroundcolor=\color{backcolour},   
    commentstyle=\color{codegreen},
    keywordstyle=\color{magenta},
    numberstyle=\tiny\color{codegray},
    stringstyle=\color{codepurple},
    basicstyle=\ttfamily\scriptsize,
    inputencoding=utf8,
    extendedchars=true,
    literate={á}{{\'a}}1 {é}{{\'e}}1 {í}{{\'i}}1 {ó}{{\'o}}1 {ú}{{\'u}}1
             {Á}{{\'A}}1 {É}{{\'E}}1 {Í}{{\'I}}1 {Ó}{{\'O}}1 {Ú}{{\'U}}1
             {ñ}{{\~n}}1 {Ñ}{{\~N}}1 {ü}{{\"u}}1 {¿}{{?`}}1 {º}{{\textordmasculine}}1
             {─}{{-}}1 {│}{{|}}1 {Σ}{{$\Sigma$}}1 {·}{{$\cdot$}}1,
    breakatwhitespace=false,
    breaklines=true,
    captionpos=b,                    
    keepspaces=true,                 
    numbers=left,                    
    numbersep=5pt,                  
    showspaces=false,                
    showstringspaces=false,
    showtabs=false,                  
    tabsize=2
}
\newcommand{\R}{\mathbb{R}}
\newcommand{\Z}{\mathbb{Z}}
\newcommand{\B}{\{ 0, 1 \}}

\newcommand{\st}{\text{s.t.}}

\newcommand{\bx}{\bm{x}}
\newcommand{\bw}{\bm{w}}
\newcommand{\bc}{\bm{c}}

\newcommand{\valueF}{\textsf{f}}
\newcommand{\transitionF}{\textsf{g}}
\newcommand{\rewardF}{\textsf{r}}

\newcommand{\horizon}{T}
\newcommand{\stages}{\mathcal{\horizon}}
\newcommand{\stagesMinusT}{\stages\setminus\{\horizon\}}
\newcommand{\feasibleSet}{\mathcal{X}}

\newcommand{\stateSpace}{\mathcal{S}}
\newcommand{\state}{s}
\newcommand{\mergestate}{\Tilde{\state}}
\newcommand{\initialState}{\state^I}
\newcommand{\dummyState}{\state^D}

\newcommand{\dd}{\mathcal{D}}
\newcommand{\nodes}{\mathcal{N}}
\newcommand{\arcs}{\mathcal{A}}

\newcommand{\maxwidth}{W}
\newcommand{\width}{w}
\newcommand{\incoming}{\arcs^\textnormal{in}}
\newcommand{\outgoing}{\arcs^\textnormal{out}}
\newcommand{\parent}{\textsf{p}}	    
\newcommand{\child}{\textsf{c}}	        

\newcommand{\rootnode}{\mathsf{r}}
\newcommand{\terminalnode}{\mathsf{t}}
\newcommand{\rt}{$\rootnode-\terminalnode$}

\newcommand{\nodeState}{\textsf{s}}

\newcommand{\arcValue}{\textsf{v}}
\newcommand{\arcLength}{\ell}

\newcommand{\merge}{\oplus}
\newcommand{\bigmerge}{\bigoplus}
\newcommand{\priorityF}{\textsf{p}}

\newcommand{\algTopDownDD}{\textsf{TopDownDD}}

\newcommand{\algDDMerge}{\textsf{MergeDDNodes}}
\newcommand{\algDDMergeGroup}{\textsf{GroupMergeDDNodes}}
\newcommand{\algDDDiscard}{\textsf{DiscardDDNodes}}

\newcommand{\algReturn}{\textbf{return}}

\tikzstyle{zero arc} = [draw,dashed, line width=0.5pt,->]
\tikzstyle{one arc} = [draw,line width=0.5pt,->]
\tikzstyle{zero arc infeasible} = [zero arc,
preaction={
	draw,gray!30!white,-,
	double=gray!30!white,
	double distance=3pt,
}]
\tikzstyle{one arc infeasible} = [one arc,
preaction={
	draw,gray!30!white,-,
	double=gray!30!white,
	double distance=3pt,
}]
\tikzstyle{optimal arc} = [draw,line width=1.3pt,->]
\tikzstyle{main node} = [circle,fill=gray!50,font=\scriptsize, inner sep=2pt]
\tikzstyle{text node} = [font=\scriptsize]

\begin{document}

\title{DD-suite: A cross-platform package to build Decision Diagrams for optimization purposes
}


\author{Antonia F. Blanco         \and
        Margarita Castro       \and
        Rodrigo Toro Icarte         
}


\institute{A. Blanco \at
              Department of Industrial and Systems Engineering, Pontificia Universidad Católica de Chile
           \and
            M. Castro \at
              Department of Industrial and Systems Engineering, Pontificia Universidad Católica de Chile \\
              \email{margarita.castro@uc.cl}           
              \and
             R. Toro Icarte \at
              Department of Computer Engineering, Pontificia Universidad Católica de Chile \\
}

\date{Received: date / Accepted: date}

\maketitle

\begin{abstract}
Decision diagrams (DDs) have become a powerful tool for discrete optimization, supporting a wide range of algorithms that span cut-generation procedures, decomposition methods, and specialized branch-and-bound searches. Despite this growth, their adoption remains limited, partly because most existing DD code is tailored to a specific algorithm or application and is therefore hard to reuse. We introduce \textit{DD-suite}, a cross-platform, open-source software package for building and manipulating DDs for discrete optimization. DD-suite is available in both Python and C++ through a shared modeling interface, and lets users construct exact, restricted, and relaxed DDs for any discrete optimization problem expressed in recursive form. The package implements the DD reduction procedure, shortest-path routines for obtaining primal and dual bounds, a visualization tool, and an extensive automated test suite. Furthermore, it includes extensive documentation, a support webpage, and ready-to-use examples for four combinatorial problems. Rather than a closed solver, DD-suite is designed as an extensible building block: users can add new construction mechanisms or run custom algorithms on top of the resulting diagram, as we illustrate with a DD-based cutting plane procedure embedded in a state-of-the-art mixed-integer programming solver. Our numerical experiments show that the C++ implementation is $5$--$6$ times faster than the Python one while producing identical diagrams, and remains within a small constant factor of \texttt{ddo}, a specialized Rust framework, confirming that DD-suite combines an accessible, extensible codebase with competitive performance.
\end{abstract}

\section{Introduction} \label{sec:intro}

In recent years, decision diagrams (DDs) have become a powerful tool for solving discrete optimization problems. This graphical structure was first proposed in the mid-1900s to represent circuits and Boolean functions \cite{akers1978binary,bryant1986graph,lee1959representation}. Still, it was not until the early 2000s that researchers started to use it for optimization purposes in mathematical programming \cite{becker2005bdds,behle2007binary}. Since then, the number of applications of DDs for optimization purposes has rapidly grown, with many applications and the creation of novel algorithms based on such structures \cite{bergman2016decision,castro2022decision,van2024introduction}. 

Specifically, DDs are graphical structures that can represent the feasibility set of most discrete optimization problems. Such representation is quite amenable to mathematical programming techniques, mainly because we can obtain the convex hull of a discrete problem using a network flow model over its corresponding DD \cite{behle2007binary}. Many authors have leveraged this characteristic to develop novel optimization algorithms, such as cut-generating procedures \cite{castro2022combinatorial,davarnia2020outer,tjandraatmadja2019target}, Benders decomposition \cite{guo2021logic,lozano2022binary,salemi2023structure}, and branch-and-price \cite{morrison2016solving,riascos2024branch} approaches, to name a few. Moreover, some researchers have created novel optimization approaches based on DDs, such as specialized branch-and-bound (B\&B) algorithms entirely based on DDs \cite{bergman2016discrete,coppe2024decision,gillard2020ddo,michel2024codd}, and the column elimination procedure \cite{karahalios2023column,van2022graph}. These techniques have been employed in a wide range of applications that include classic combinatorial problems (e.g., independent set and graph coloring \cite{bergman2014optimization,van2022graph}), sequencing and vehicle routing problems \cite{cire2013multivalued,kinable2017hybrid,castro2020mPDTSP}, and healthcare applications \cite{cire2019network,guo2021logic,riascos2024branch}, to name a few. 

Despite the wide range of algorithms and applications based on DDs, their use is still quite limited, and only a few research groups have adopted such techniques. We believe this lack of diversity is mainly because DD-based algorithms require extensive coding to create and manipulate such DDs, which, until now, have been developed independently for each specific project. While several researchers have made their code publicly available, such code is usually tuned to their particular algorithm and application, which makes it very hard to modify or extend for other usages. To the best of our knowledge, two of the more general-purpose codes available for DDs in optimization are the DD-based B\&B solvers (e.g., \texttt{ddo} \cite{gillard2020ddo} and CODD \cite{michel2024codd}), where users can easily include new problems that can be solved with their solver. However, these codes are mainly developed as solvers and are hard to adapt for other purposes (see Section \ref{sec:comparison} for further details and discussion). 

Given the lack of software tools for DD-based optimization research, this work presents \textit{DD-suite}, a cross-platform open-source code\footnote{Available at \url{https://github.com/MargaritaCastro/dd-suite}.} to create and manipulate DDs for discrete optimization purposes. This software suite is available in Python and C++ and allows users to easily create three types of DDs (exact, restricted, or relaxed) for any discrete optimization problem they prefer. Our main goal is to provide an intuitive, clean, and clear code that researchers can use to build DDs and create novel algorithms. To do so, DD-suite follows several good coding practices from software engineering \cite{martin2009clean}, including clear on-code documentation and a support webpage with detailed tutorials and examples. 

DD-suite is designed to be easily adapted to different usages, so researchers can employ DDs in any form they like. As an illustrative example, we provide a simple cutting plane implementation using DD-suite, in which we implemented existing DD-based cut-generation procedures (see Section \ref{sec:cutting-planes} for further details). Moreover, our implementation allows users to test and explore new ways to construct DDs by changing key functions of their problem specifications or extending DD-suite to incorporate additional construction mechanisms. Furthermore, we include extensive automatic testing so users can safely work with DD-suite without fearing breaking the existing code. Lastly, DD-suite includes examples for four different combinatorial problems (i.e., knapsack, independent set, set cover, and sequencing), together with a second-order cone (SOC) knapsack example used in our cutting plane experiments, all with intuitive main files for users to play with.

Our numerical experiments confirm that this flexibility does not come at the expense of performance: the C++ and Python implementations produce identical diagrams, with C++ being $5$--$6$ times faster and remaining within a small constant factor of \texttt{ddo}, a state-of-the-art Rust framework heavily optimized for speed. We further use these experiments to showcase the extensibility of DD-suite (i.e., by adding an alternative DD construction mechanism) and its research value (i.e., by implementing DD-based cutting planes that improve a state-of-the-art commercial solver).

\paragraph{Contributions.}
Overall, the main contributions of this work are:
\begin{itemize}
    \item A clean, intuitive, and easily extended Python and C++ code to create and manipulate DDs for optimization purposes. Our implementation supports exact, relaxed, and restricted DDs and includes the DD reduction procedure, a shortest-path implementation to obtain optimality bounds, a DD visualization tool,  and many other functionalities to get information about the DD and the supported algorithms. 
    \item On-code documentation and a support webpage with tutorials and useful information on DD-suite usage. 
    \item Extensive automatic testing to safely adapt and modify DD-suite for different usages without jeopardizing existing functionalities. 
    \item A cutting plane example to illustrate how to use DD-suite for advanced optimization purposes and how to integrate it with mathematical programming solvers (e.g., Gurobi). 
    \item Experiments comparing the efficiency of both code alternatives (Python vs. C++) and existing DD-based solvers. 
\end{itemize}

\paragraph{Paper Structure.}
The rest of the paper is organized as follows. Section \ref{sec:comparison} reviews the publicly available DD software most closely related to DD-suite and positions our work with respect to them. Section \ref{sec:background} provides a background on DDs for discrete optimization, where we describe the main algorithms implemented in DD-suite. Section \ref{sec:software} details the basic features of the software (i.e., creating and manipulating DDs) and provides documentation and webpage details. Section \ref{sec:developers} covers advanced usage aimed at developers, including the automatic testing framework, the cutting plane example, and how to extend DD-suite for new problems and construction mechanisms. Section \ref{sec:experiments} presents our numerical experiments, comparing the efficiency of the Python and C++ implementations against each other and against existing DD-based solvers, and illustrates the use of DD-suite for cutting planes and for prototyping new relaxation mechanisms. Finally, Section \ref{sec:conclusions} concludes the paper and outlines future work directions.

\section{Literature Review} \label{sec:comparison}

We now review available tools that are closest in spirit to DD-suite, namely those that provide a somewhat general interface for modeling problems and building DDs. We focus on the two general-purpose DD-based B\&B solvers, \texttt{ddo} \cite{gillard2020ddo} and CODD \cite{michel2024codd}, and on the DD-compilation framework Haddock \cite{gentzel2020haddock}. Table \ref{tab:comparison} summarizes the comparison, where a check mark (\checkmark) denotes a supported feature and a dash (--) an unsupported one.

\begin{table}[b]
\centering
\caption{DD-suite versus the most closely related software.}
\label{tab:comparison}
\setlength{\tabcolsep}{5pt}
\renewcommand{\arraystretch}{1.15}
\begin{tabular}{lcccc}
\toprule
 & \texttt{ddo} & CODD & Haddock & DD-suite \\
\midrule
Language(s)              & Rust          & C++           & C++ (MiniCP)   & Python \& C++ \\
Primary role             & B\&B solver     & B\&B solver     & DD compilation / CP & DD library \\
Problem specification    & DP recursion  & DP (lambdas)  & Transition system & DP recursion \\
Exact / restricted / relaxed DDs & \checkmark & \checkmark & \checkmark & \checkmark \\
Primal \& dual bounds    & \checkmark    & \checkmark    & \checkmark      & \checkmark \\
Parallel search          & \checkmark    & --            & --              & -- \\
Exposes DD for reuse     & --            & --            & --              & \checkmark \\
Cutting planes / MP integration & --     & --            & --              & \checkmark \\
Cross-language interface & --            & --            & --              & \checkmark \\
Documentation \& tutorials & README      & Limited       & Limited         & Extensive \\
Automatic test suite     & \checkmark    & --            & --              & \checkmark \\
Open source              & \checkmark    & \checkmark    & \checkmark      & \checkmark \\
\bottomrule
\end{tabular}
\end{table}

The \texttt{ddo} framework \cite{gillard2020ddo} is a generic and efficient Rust library for optimization based on decision diagrams. It implements the DD-based B\&B algorithm of \cite{bergman2016discrete}, where relaxed and restricted DDs provide the dual and primal bounds that drive the search, and it supports parallel B\&B. Subsequent works have improved its performance through better filtering rules \cite{gillard2021improving} and a caching mechanism for dominance and suboptimality detection \cite{coppe2024decision}. The framework includes several examples (e.g., knapsack, independent set, MaxSAT, and MaxCut) and is heavily optimized for raw performance, which makes it a natural performance baseline for our experiments (see Section \ref{sec:exp-efficiency}); on the other hand, its design centers on the B\&B solver rather than on exposing the DDs for general use.

CODD \cite{michel2024codd} is a more recent C++ DD-based solver. Like \texttt{ddo}, it performs a DD-based B\&B search, but it lets users specify the dynamic programming (DP) model very concisely through C++ lambda functions, sitting conceptually between a general DP solver, such as the domain-independent dynamic programming (DIDP) framework of \cite{kuroiwa2023domain}, and a specialized DD solver such as \texttt{ddo}. As with \texttt{ddo}, CODD is primarily designed to be used as a solver, and its DD components are not meant to be repurposed for other DD-based techniques.

Haddock \cite{gentzel2020haddock} is a language and architecture for DD compilation built on top of the MiniCP system (C++). Instead of an explicit DP recursion, problems are described through a labeled transition system, from which DDs are compiled. Haddock is mainly focused on DD propagation for constraint programming (CP), with later work adding the computation of optimization bounds \cite{gentzel2023optimization}. Its scope is therefore different from ours: it targets DD-based propagation within a CP solver rather than providing a general-purpose toolkit to build and manipulate DDs for mathematical-programming-style algorithms.

Beyond these general frameworks, DDs have been embedded in several other optimization and CP tools. For instance, the peel-and-bound algorithm \cite{rudich2022peel} builds relaxed DDs through \emph{separation} (i.e., starting from a small, weak diagram and iteratively refining it by splitting nodes) to generate strong dual bounds, an alternative construction paradigm to the top-down procedure used by the solvers above and by DD-suite. Performance has also been pushed along a different axis: \cite{tardivo2026complete} accelerates the construction of DDs for DP models on GPUs, expanding and filtering states in parallel to obtain substantial speedups on hard sequencing problems. In the CP community, DDs are a standard data structure, and efficient libraries of DD operations (e.g., reduction, intersection, and union) have been developed to build CP models \cite{perez2015efficient}. These tools, however, are again specialized to their target use (dual-bound computation, GPU-accelerated search, or CP propagation) rather than offering a general toolkit to build and manipulate DDs for arbitrary optimization algorithms.

A common feature of all these tools, and of DD-suite, is that the user specifies a problem through a recursive, DP model (i.e., a state representation together with transition, merging, and discarding rules), from which the DD is automatically constructed. The main difference lies in their intended use: \texttt{ddo} and CODD are \emph{end-to-end solvers} in which the DD machinery is internal to a B\&B search, whereas DD-suite is a \emph{library} that exposes the DD construction and the resulting graph so that they can be reused as a building block for arbitrary DD-based algorithms.


\section{Background on Decision Diagrams} \label{sec:background}

This section provides a brief overview of DDs for optimization. We refer the reader to \cite{bergman2016decision,castro2022decision,van2024introduction} for further details on DDs and their usage.  In what follows, we use calligraphic font for sets, uppercase letters for the sets' cardinality, and lowercase letters for the sets' elements if possible (e.g., $a\in \arcs$ is an element of set $\arcs$ with cardinality $|\arcs|=A$). We also employ bold letters for vectors and sans serif font for functions (e.g., $\bx=(x_1,x_2) \in \Z^2$ and $\textsf{a}(\cdot): \Z \rightarrow \R$ is a function). 

Given the strong relationship between DP and DDs \cite{hooker2013decision}, we consider recursive formulations of discrete problems to explain DDs and as the problem specification for our software (see Section \ref{sec:software} for more details). Specifically, we consider discrete problems that can be represented sequentially and have a finite number of variables with bounded domain (i.e., $\bx=(x_1,...,x_\horizon) \in \feasibleSet \subset \Z^\horizon$, where  $\horizon \in \Z$ and $\feasibleSet$ is bounded).

Consider $\stages=\{1,...,\horizon\}$ as the set of stages for the recursive formulation, $\initialState$ as the initial state of the system, and $\stateSpace$ as the set of reachable states starting from $\initialState$.  Abusing the notation, we use  $\stateSpace_t$ to represent reachable states at stage $t \in \stages$ and $\feasibleSet_t(s)$ as the set of feasible assignments of $x_t$ at state $s\in \stateSpace_t$. Then, the Bellman equations describing a maximization problem for all $s \in \stateSpace$ are:
\begin{equation}\label{eq:bellman}
    \valueF_t(\state) = \begin{cases}
    \max_{x_t \in \feasibleSet_t(s)} \{ \rewardF_t (\state, x_t) + \valueF_{t+1}( \transitionF_t(\state, x_t) ) \}, &  s\in \stateSpace_t, \; t\in \stagesMinusT, \\
    \max_{x_t \in \feasibleSet_t(s)} \{ \rewardF_t (\state, x_t) \}, &  s\in \stateSpace_T,\;  t = \horizon.
    \end{cases}
\end{equation}
The value function at stage $t \in \stages$ is given by $\valueF_t(\cdot): \stateSpace \rightarrow \R$, while $\rewardF_t(\cdot, \cdot): \stateSpace_t \times \feasibleSet_t \rightarrow \R$ is the immediate reward of choosing a feasible value for variable $x_t \in \feasibleSet_t(s)$ at state $\state \in \stateSpace_t$. Lastly, $\transitionF_t(\cdot, \cdot): \stateSpace_t \times \feasibleSet_t \rightarrow \stateSpace_{t+1}$ corresponds to the transition function of choosing a value for $x_t$ at state $s$. Then, the optimal solution for this recursive model is given by $\valueF_1(\initialState)$.

\begin{example}\label{example:recursiveKnapsack}
To better explain the main concepts related to DDs, we borrow the knapsack example in \cite{castro2022decision} as our running example throughout the paper. Consider a knapsack problem with four variables, weight vector $\bw = (7,5,4,1)$, and reward vector $\bc=(4,2,5,1)$, that is,
\begin{align*}
 \max\; & 4x_1 + 2x_2 + 5x_3 + x_4 \\
 \st\; & 7x_1 + 5x_2 + 4x_3 +x_4 \leq 8, \\
 & \bx=(x_1,...,x_4) \in \B^4.
\end{align*}

A recursive formulation of the problem considers choosing one item per stage in lexicographic order. Each state represents the current load of the knapsack (i.e., an integer value), and the initial state is given by $\initialState = 0$. For a given state $s \in \stateSpace_t$ and stage $t \in \stages$, the set of feasible values of $x_t$ is $\feasibleSet_t(s)= \{x\in \B:\; s + w_t x \leq 8 \}$. Then, the transition function for $s \in \stateSpace_t$ and $x_t \in \feasibleSet_t(s)$ is  $\transitionF_t(s,x_t) = s + w_tx_t$. Note that the immediate reward function is independent of the state and, thus, is given by $\rewardF_t(\state, x_t) = c_tx_t$. Therefore, the Bellman equations governing this problem are:
\begin{equation*}
    \valueF_t(\state) = \begin{cases}
    \max_{x_t \in \feasibleSet_t(s)} \{ c_t x_t + \valueF_{t+1}( \state + w_t x_t ) \}, &  s\in \stateSpace_t, \; t\in \{1,...,3\}, \\
    \max_{x_t \in \feasibleSet_t(s)} \{ c_t x_t \}, &  s\in \stateSpace_4,\;  t = 4,
    \end{cases}
    \quad \forall s \in \stateSpace.
\end{equation*}
\hfill $\square$
\end{example}

\subsection{Exact Decision Diagrams}

DDs are graphical representations of the feasible set of a discrete problem. Given a recursive formulation, a DD can be seen as a compact representation of the transition graph of such a model, where nodes map to states and arcs map to feasible variable assignments. We now describe the main components and properties of such a graphical structure. 

Consider the recursive model given by equations \eqref{eq:bellman}. Formally, a DD is a directed acyclic layered graph $\dd=(\nodes,\arcs)$, where $\nodes$ and $\arcs$ are the set of nodes and arcs, respectively. The set of nodes is partitioned into $\horizon+1$ layers $\nodes = \{\nodes_1,..., \nodes_{\horizon+1} \}$ such that each node $n \in \nodes_t$ is associated with a single state $s\in \stateSpace_t$, for all $t \in \stages$. We use the function $\nodeState(\cdot):\nodes \rightarrow \stateSpace$ to represent the state associated with a node. In particular, the first node layer has a single node, called the \textit{root node} (i.e.,  $\nodes_1=\{ \rootnode \}$), which is associated with the initial state (i.e., $\nodeState(\rootnode) = \initialState$). The last node layer also has a single node called the \textit{terminal node} (i.e., $\nodes_{\horizon + 1} = \{\terminalnode\}$) and represents all the terminal states for the recursive model. For simplicity, we associate a dummy state $\dummyState$ to the terminal node (i.e., $\nodeState(\terminalnode) = \dummyState$). 

The set of arcs is partitioned into $\horizon$ layers, $\arcs= \{\arcs_1,...,\arcs_T\}$, such that an arc $a \in \arcs_t$ emanates from its parent node $\parent_a \in \nodes_{t}$ and points to its child node $\child_a \in \nodes_{t+1}$, for all $t \in \stages$. Each arc $a \in \arcs_t$, for $t\in \stages$, has a value $\arcValue_a \in \feasibleSet_t( \nodeState(\parent_a))$ corresponding to a feasible assignment of $x_t$; and its child node is associated with the state given by the transition function (i.e., $\nodeState(\child_a) = \transitionF_t(\nodeState(\parent_a), \arcValue_a)$). Given that the last node layer has a single node, the child node for all arcs in the last layer is the terminal node (i.e., $\child_a = \terminalnode$ for all $a\in \arcs_\horizon$), and thus, the transition function should return the dummy state $\dummyState$ in such cases. Abusing the notation, we use $\outgoing(n)$ and $\incoming(n)$ to represent the set of outgoing and incoming arcs for any node $n \in \nodes$ (i.e., $a \in \outgoing(n)$ and $a \in \incoming(n')$ if and only if $\parent_a = n$ and $\child_a = n'$).

As illustrated in Example \ref{example:exactDD}, a DD is indeed a compact representation of the state transition graph of a recursive model. Moreover, we say that a DD is \textit{exact} if a one-to-one mapping exists between the feasible assignments of $\bx$ and the \rt\ paths on the DD (e.g., a DD following the description provided in this section). In such case, we can obtain the optimal solution of the recursive model by solving the longest-path problem from $\rootnode$ to $\terminalnode$  where the length of each arc $a \in \arcs_t$, $t \in \stages$, is given by the reward function as $\arcLength_a = \rewardF_t(\nodeState(\parent_a), \arcValue_a)$.

One of the main advantages of DDs is that we can further compress the graphical structure without sacrificing any feasible assignment. This is done using a reduction procedure \cite{bryant1992symbolic} that merges isomorphic subgraphs by traversing each node layer once, starting from $\nodes_{\horizon+1}$, and runs in polynomial time with respect to the DD size (i.e., number of nodes and arcs). The procedure returns the smallest exact DD at the cost of some nodes representing two or more states. Moreover, we can use any longest-path algorithm on the reduced graph to obtain the optimal solution if the reward function is not state-dependent (e.g., linear).

\begin{figure}
    \centering
    \begin{tikzpicture}[->,>=stealth',shorten >=1pt,auto,node distance=1cm,
thick]        
\node[text node] (r) at (0,0) {$\quad$};
\node[text node] (t) at (0,-4.1) {$\quad$};

\node[text node] (l1) at (0,-0.5) {$x_1$:};
\node[text node] (l2) at (0,-1.5) {$x_2$:};
\node[text node] (l4) at (0,-2.5) {$x_3$:};
\node[text node] (l5) at (0,-3.5) {$x_4$:};
\end{tikzpicture} \hfill
    \begin{tikzpicture}[->,>=stealth',shorten >=1pt,auto,node distance=1cm,
thick]        
\node[main node] (r) at (0,0) {0};
\node[main node] (u1)  at (-1.2,-1)  {0};
\node[main node] (u2)  at (1.2,-1)  {7};
\node[main node] (u3)  at (-1.2,-2) {0};
\node[main node] (u4)  at (0,-2) {5};
\node[main node] (u5)  at (1.2,-2) {7};
\node[main node] (u6)  at (-1.5,-3) {0};
\node[main node] (u7)  at (-0.5,-3) {4};
\node[main node] (u8)  at (0.5,-3) {5};
\node[main node] (u9)  at (1.5,-3) {7};
\node[main node] (u10) at (-1.5,-4) {0};
\node[main node] (u11) at (-1,-4) {1};
\node[main node] (u12) at (-0.5,-4) {4};
\node[main node] (u13) at (0,-4) {5};
\node[main node] (u14) at (0.5,-4) {6};
\node[main node] (u15) at (1,-4) {7};
\node[main node] (u16) at (1.5,-4) {8};


\path[every node/.style={font=\sffamily\small}]
(r) 
edge[zero arc, optimal arc] node [left] {} (u1)
edge[one arc] node [left] {} (u2)
(u1) 
edge[zero arc, optimal arc] node [right] {} (u3)
edge[one arc] node [right] {} (u4)
(u2)
edge[zero arc] node [left] {} (u5)
(u3)
edge[zero arc] node [right] {} (u6)
edge[one arc, optimal arc] node [right] {} (u7)
(u4)
edge[zero arc] node [left] {} (u8)
(u5)
edge[zero arc] node [left] {} (u9)
(u6)
edge[zero arc] node [left] {} (u10)
edge[one arc] node [left] {} (u11)
(u7)
edge[zero arc] node [left] {} (u12)
edge[one arc, optimal arc] node [left] {} (u13)
(u8)
edge[zero arc] node [left] {} (u13)
edge[one arc] node [left] {} (u14)
(u9)
edge[zero arc] node [left] {} (u15)
edge[one arc] node [left] {} (u16);

\end{tikzpicture} \hfill
    \begin{tikzpicture}[->,>=stealth',shorten >=1pt,auto,node distance=1cm,
thick]        
\node[main node] (r) at (0,0) {$\;\rootnode\;$};
\node[main node] (u1)  at (-1,-1)  {$0$};
\node[main node] (u2)  at (1,-1)  {$7$};
\node[main node] (u3)  at (-1,-2) {$0$};
\node[main node] (u4)  at (0,-2) {$5$};
\node[main node] (u5)  at (1,-2) {$7$};
\node[main node] (u6)  at (-1.2,-3) {$0$};
\node[main node] (u7)  at (-0.3,-3) {$4$};
\node[main node] (u8)  at (0.3,-3) {$5$};
\node[main node] (u9)  at (1.2,-3) {$7$};
\node[main node] (t) at (0,-4) {$\;\terminalnode\;$};

\path[every node/.style={font=\sffamily\small}]
(r) 
edge[zero arc, optimal arc] node [left] {} (u1)
edge[one arc] node [left] {} (u2)
(u1) 
edge[zero arc, optimal arc] node [right] {} (u3)
edge[one arc] node [right] {} (u4)
(u2)
edge[zero arc] node [left] {} (u5)
(u3)
edge[zero arc] node [right] {} (u6)
edge[one arc, optimal arc] node [right] {} (u7)
(u4)
edge[zero arc] node [left] {} (u8)
(u5)
edge[zero arc] node [left] {} (u9)
(u6)
edge[zero arc,  bend right=20] node [left] {} (t)
edge[one arc,  bend left=20] node [left] {} (t)
(u7)
edge[zero arc,  bend right=20] node [left] {} (t)
edge[one arc,  bend left=20, optimal arc] node [left] {} (t)
(u8)
edge[zero arc,  bend right=20] node [left] {} (t)
edge[one arc,  bend left=20] node [left] {} (t)
(u9)
edge[zero arc,  bend right=20] node [left] {} (t)
edge[one arc,  bend left=20] node [left] {} (t);

\end{tikzpicture} \hfill
    \begin{tikzpicture}[->,>=stealth',shorten >=1pt,auto,node distance=1cm,
thick]        
\node[main node] (r) at (0,0) {$\;\rootnode\;$};
\node[main node] (u1)  at (-1,-1)  {$\quad$};
\node[main node] (u2)  at (1,-1)  {$\quad$};
\node[main node] (u3)  at (-1,-2) {$\quad$};
\node[main node] (u4)  at (1,-2) {$\quad$};
\node[main node] (u5)  at (0,-3) {$\quad$};
\node[main node] (t) at (0,-4) {$\;\terminalnode\;$};

\path[every node/.style={font=\scriptsize}]
(r) 
edge[zero arc, optimal arc] node [left] {0} (u1)
edge[one arc] node [right] {4} (u2)
(u1) 
edge[zero arc, optimal arc] node [left] {0} (u3)
edge[one arc] node [right] {2} (u4)
(u2)
edge[zero arc] node [right] {0} (u4)
(u3)
edge[zero arc, bend right=20] node [left] {0} (u5)
edge[one arc, bend left=20, optimal arc] node [right] {5} (u5)
(u4)
edge[zero arc] node [right] {0} (u5)
(u5)
edge[zero arc,  bend right=45] node [left] {0} (t)
edge[one arc,  bend left=45, optimal arc] node [right] {1} (t);

\end{tikzpicture}
    \caption{State transition graph, exact DD, and its reduced version.}
    \label{fig:exactDDKnapsack}
\end{figure}
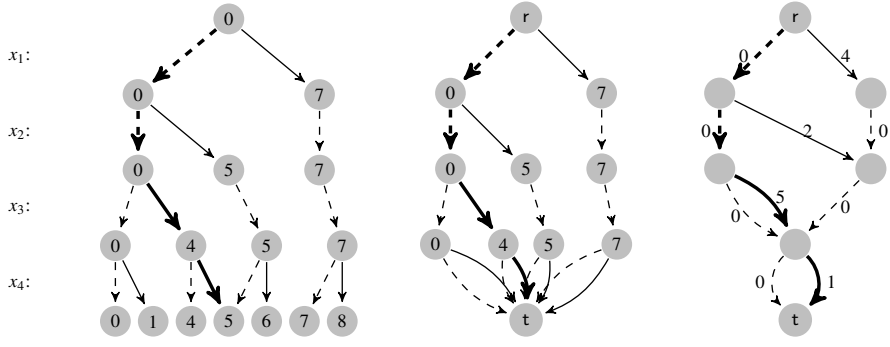

\begin{example}\label{example:exactDD}
    Following the knapsack problem described in Example \ref{example:recursiveKnapsack}, Figure \ref{fig:exactDDKnapsack} presents the state transition graph for the model (left), an exact DD (middle), and its reduced version (right). In all graphs, the numbers inside the circles correspond to states, and arcs represent variable assignments (solid arcs for $x_t = 1$ and dashed arcs for $x_t = 0$). Also, numbers next to arcs are their associated length given by the reward function, and the bold paths correspond to the optimal solution to the problem.
    
    We observe that the middle graph is an exact DD where each node is associated with a single state. The main difference with the transition graph is that there is a single node in the last layer (i.e., the terminal node $\terminalnode$). The reduced DD is obtained by applying the reduction procedure \cite{bryant1992symbolic} to the middle DD. This rightmost DD has considerably fewer nodes (i.e., at most two per layer), and it preserves all the paths (i.e., feasible assignments) present in the transition graph and the middle DD. Lastly, we can obtain the optimal solution from all graphical structures using a longest-path algorithm where $\arcLength_a = c_t\arcValue_a$ for all $a \in \arcs_t$, $t\in \stages$, which corresponds to $\bx^*=(0,0,1,1)$ with value $6$.
\hfill $\square$
\end{example}

\subsection{Approximated Decision Diagrams}

One of the main drawbacks of exact DDs is that they suffer from the \textit{curse of dimensionality} \cite{bellman1957dynamic}; that is, the number of nodes and arcs grows exponentially with the number of variables for most problems. While the reduction procedure can somewhat mitigate this effect, its final size can be exponential in the number of variables. To overcome this issue, researchers have proposed two types of DDs: underapproximations (i.e., restricted DDs) or overapproximations (i.e., relaxed DDs) of the feasibility set of the original problem. We now describe each approximated DD and provide general insights.

\subsubsection{Restricted DDs}

Restricted DDs were first introduced as a general procedure to obtain high-quality solutions for discrete optimization problems \cite{bergman2014heuristic}, but they can also be used, for instance, to create relaxations for bilevel problems \cite{vasquez2025bilevel}. The main idea is to construct a limited-size DD by discarding nodes during construction that could lead to poor-quality solutions. The most common approach to limit a DD size is through its maximum allowed width, that is, the maximum number of nodes per layer. Formally, given a DD $\dd=(\nodes,\arcs)$, we define the width of a node layer $\nodes_t$ as the number of nodes for such layer, that is, $\width(\nodes_t) = |\nodes_t|$ for all $t\in \stages$. Then, the maximum width of a DD is given by $\maxwidth= \max_{t \in \stages} \{\width(\nodes_t)\}$. 

For a given maximum width $\maxwidth$, we build a restricted DD following a procedure similar to the one used for an exact DD (i.e., one node per state in each stage). If the number of nodes in a layer is larger than $\maxwidth$, then we choose $\maxwidth$ nodes that might lead to an optimal solution and discard the rest. This discarding procedure is usually heuristic and problem-specific (see, for instance, Example \ref{example:approximatedDD}). The resulting restricted DD has a limited size where the number of nodes per layer is at most $\maxwidth$, and its paths represent a subset of all possible solutions to the original problem. Therefore, we can optimize over the resulting restricted DD and obtain a primal bound of the original problem. More details of the construction mechanism can be found in Section \ref{sec:dd-construction}.

\subsubsection{Relaxed DDs}  \label{sec:relaxedDDs}

Relaxed DDs are graphical structures that represent a discrete relaxation of an optimization problem and, thus, can be used to obtain dual bounds \cite{andersen2007constraint}. Similar to restricted DDs, relaxed DDs can be constructed by limiting the number of nodes in each layer (i.e., its maximum width). However, instead of discarding nodes when the size is exceeded, these DDs merge nodes such that a node represents two or more states. By doing so, the resulting DD might incorporate paths that represent infeasible assignments, thus overapproximating the original problem. 

While the main idea of relaxed DDs is quite simple (i.e., merging nodes and states), it is crucial to define the merging mechanism and its implications on the recursive model (i.e., states, transition function, and reward function) to guarantee that the resulting DD indeed relaxes the original problem. Hooker \cite{hooker2017job} formalizes the key properties needed for a merging operator to obtain a relaxed DD. Specifically, given two states $\state, \state' \in \stateSpace_t$ for any $t\in \stages$, we define a \textit{merged state} $\mergestate = \state\merge\state'$, where $\merge$ is a merge operator. This merge operator defines a proper relaxation if the following conditions are satisfied:
\begin{enumerate}[label=(\roman*)]
    \item The feasibility set of $\mergestate$ contains all the feasible assignments associated with  $\state$ and $\state'$, that is, $\feasibleSet_t(\state)\cup \feasibleSet_t(\state') \subseteq \feasibleSet_t(\mergestate)$. \label{cond:feasibleset}
    \item The immediate reward at state $\mergestate$ is greater than or equal to the immediate reward at $\state$ and $\state'$ for any feasible assignment, that is, given $x\in \feasibleSet_t(\state)$ and $x' \in \feasibleSet_t(\state')$ we have $\rewardF_t(\state,x) \leq \rewardF_t(\mergestate,x)$ and $\rewardF_t(\state', x')\leq \rewardF_t(\mergestate, x')$. \label{cond:reward}
    \item If $\mergestate$ \textit{relaxes} state $\state$ (i.e., $\mergestate$ satisfies \ref{cond:feasibleset} and \ref{cond:reward}), then $\transitionF_t(\mergestate, x)$ relaxes state $\transitionF_t(\state, x)$ for all $x\in \feasibleSet_t(\state)$. \label{cond:transition}
\end{enumerate}

As shown in \cite{hooker2017job}, these three conditions are enough to guarantee that the merging operator $\merge$ will provide a relaxed DD that indeed overapproximates the feasible set of the problem and, thus, can provide dual bounds. As in the restricted DD case, the merging operator is problem-specific, and choosing which nodes to merge is usually done heuristically (see Example \ref{example:approximatedDD} for its application to a knapsack instance). 

\begin{figure}
    \centering
    \begin{tikzpicture}[->,>=stealth',shorten >=1pt,auto,node distance=1cm,
thick]        
\node[text node] (r) at (0,0) {$\quad$};
\node[text node] (t) at (0,-4.1) {$\quad$};

\node[text node] (l1) at (0,-0.5) {$x_1$:};
\node[text node] (l2) at (0,-1.5) {$x_2$:};
\node[text node] (l4) at (0,-2.5) {$x_3$:};
\node[text node] (l5) at (0,-3.5) {$x_4$:};
\end{tikzpicture} \hspace{2em}
    \begin{tikzpicture}[->,>=stealth',shorten >=1pt,auto,node distance=1cm,
thick]        
\node[main node] (r) at (0,0) {$\;\rootnode\;$};
\node[main node] (u1)  at (-1,-1)  {$0$};
\node[main node] (u11)  at (1,-1)  {$7$};
\node[main node] (u2)  at (-1,-2) {$5$};
\node[main node] (u3)  at (1,-2) {$7$};
\node[main node] (u5)  at (-1,-3) {$5$};
\node[main node] (u6)  at (1,-3) {$7$};
\node[main node] (t) at (0,-4) {$\;\terminalnode\;$};

\path[every node/.style={font=\sffamily\small}]
(r) 
edge[zero arc] node [left] {} (u1)
edge[one arc, optimal arc] node [left] {} (u11)
(u1) 
edge[one arc] node [right] {} (u2)
(u11)
edge[zero arc, optimal arc] node [right] {} (u3)
(u2)
edge[zero arc] node [right] {} (u5)
(u3)
edge[zero arc, optimal arc] node [left] {} (u6)
(u5)
edge[zero arc,  bend right=20] node [left] {} (t)
edge[one arc,  bend left=20] node [left] {} (t)
(u6)
edge[zero arc,  bend right=20] node [left] {} (t)
edge[one arc,  bend left=20, optimal arc] node [left] {} (t);

\end{tikzpicture} \hspace{5em}
    \begin{tikzpicture}[->,>=stealth',shorten >=1pt,auto,node distance=1cm,
thick]
\node[main node] (r) at (0,0) {$\;\rootnode\;$};
\node[main node] (u1)  at (-1,-1)  {$0$};
\node[main node] (u2)  at (1,-1)  {$7$};
\node[main node] (u3)  at (-1,-2)  {$0$};
\node[main node] (u4)  at (1,-2)  {$7$};
\node[main node] (u5)  at (-1,-3) {$0$};
\node[main node] (u6)  at (1,-3) {$7$};
\node[main node] (t) at (0,-4) {$\;\terminalnode\;$};

\path[every node/.style={font=\sffamily\small}]
(r)
edge[zero arc infeasible, optimal arc] node [left] {} (u1)
edge[one arc] node [left] {} (u2)
(u1)
edge[zero arc, bend right=30] node [right] {} (u3)
edge[one arc infeasible, optimal arc, bend left=30] node [right] {} (u3)
(u2)
edge[zero arc] node [left] {} (u4)
(u3)
edge[zero arc, bend right=30] node [left] {} (u5)
edge[one arc infeasible, optimal arc, bend left=30] node [right] {} (u5)
(u4)
edge[zero arc] node [left] {} (u6)
(u6)
edge[zero arc, bend right=20] node [left] {} (t)
edge[one arc, bend left=20] node [left] {} (t)
(u5)
edge[zero arc, bend right=20] node [left] {} (t)
edge[one arc infeasible, optimal arc, bend left=20] node [left] {} (t);

\end{tikzpicture}
    \caption{Restricted and relaxed DDs for the knapsack example.}
    \label{fig:approxDDKnapsack}
\end{figure}
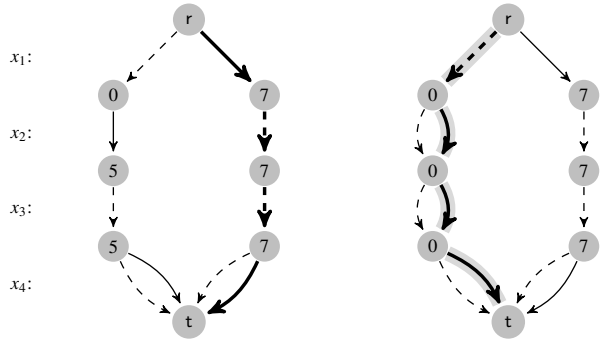

\begin{example}\label{example:approximatedDD}
    Following the knapsack instance described in the previous examples, Figure \ref{fig:approxDDKnapsack} depicts a restricted DD (left) and a relaxed DD (right) of maximum width $\maxwidth = 2$ for such a problem. The restricted DD uses a discarding heuristic where we keep the nodes with the highest knapsack load. Indeed, the third layer has only two nodes and discards the node with $\state=0$.  Note that every path in the restricted DD corresponds to a feasible assignment of the problem, but some solutions (e.g., $\bx=(0,0,0,0)$) are not represented in the DD. The optimal solution (i.e., a primal bound for the original problem) is depicted by the bold path in the figure, which represents solution $\bx=(1,0,0,1)$ with value $5$.

    To define a relaxed DD, the merging operator keeps the smallest load among the merged nodes, that is, 
    $\mergestate = \state \merge\state' = \min\{\state, \state'\}$ for any two states $\state, \state' \in \stateSpace_t$ (for some $t\in \stages$). Since a smaller load leaves more room in the knapsack, the merged state contains all the feasible assignments of $\state$ and $\state'$; and, since the reward function is state-independent, we keep it as it is. Note that tracking only the smallest load suffices to obtain a valid relaxation precisely because the reward does not depend on the state. This merging operator follows conditions \ref{cond:feasibleset}-\ref{cond:transition} and, thus, provides a relaxed DD.

    The relaxed DD in Figure \ref{fig:approxDDKnapsack} utilizes a merging heuristic (which nodes to merge) that merges nodes with the lowest knapsack load. Note that all feasible solutions are represented in the graph, but some paths (i.e., the shaded path) are infeasible. Indeed, the optimal solution corresponds to the shaded path associated with $\bx=(0,1,1,1)$ and value $8$. Therefore, the restricted DD provides a primal bound of 5 and the relaxed DD a dual bound of 8 (recall that the optimal value is 6).   \hfill $\square$
\end{example}

The main advantage of these two DD approximations is that the user can easily control the quality of the bounds by simply adjusting the maximum width or trying different discarding/merging strategies \cite{bergman2014optimization}. Moreover, they can be used to find the optimal solution to a problem using a searching procedure such as B\&B \cite{bergman2016discrete,coppe2024decision} or $A^*$ search \cite{castro2019relaxedbdds,castro2020SolvingDF},  or be used as a complement to existing mathematical programming solvers (e.g., to generate cuts \cite{castro2022combinatorial,davarnia2020outer,tjandraatmadja2019target} or solve the subproblems of a Lagrangian decomposition approach \cite{bergman2015lagrangian,castro2020mPDTSP}, among many others). We refer the reader to \cite{castro2022decision,van2024introduction} for further recent examples of how to employ exact and approximated DDs to solve optimization problems.

\subsection{Construction Mechanism}\label{sec:dd-construction}

We now formally present the top-down construction mechanism for exact and approximated DDs. We note that DD-suite uses this procedure to construct DDs, which can be easily tailored to different problems, as we discuss in Section \ref{sec:software}. As previously mentioned, the procedure employs a recursive formulation for a discrete problem with a fixed number of variables (i.e., given by the Bellman equations \eqref{eq:bellman}) and constructs all possible states at each stage.

Specifically, Algorithm \ref{alg:dd_topdown} describes the complete top-down construction procedure. The algorithm first initializes the graphical structure and sets the states for the root and terminal node accordingly (lines 2--3). We then iterate over each node layer in lexicographic order, and for each node, we create an arc for every feasible value associated with that node's state (lines 4--6). We then create the nodes in the next layer using the recursive model's transition function, ensuring that two nodes cannot be associated with the same state (i.e., lines 7--11). If the maximum width of the following layer is exceeded, we can merge or discard nodes depending on the type of DD we want (lines 12--13). Lastly, we create the arcs in the last layer and direct them to the terminal node. 

\begin{algorithm}[tb]
	\caption{DD Top-Down Construction} \label{alg:dd_topdown}
	\begin{algorithmic}[1]
		\Procedure{\algTopDownDD}{$\maxwidth$}
		\State Create DD $\dd=(\nodes, \arcs)$ with $\horizon + 1$ empty node layers
		\State Create the root and terminal node, i.e.,  $\rootnode\in \nodes_{1}$ with $\nodeState(\rootnode) = \initialState$ and $\terminalnode \in \nodes_{\horizon+1}$ with $\nodeState(\terminalnode) = \dummyState$
		\For{ $t \in \stagesMinusT$}
		\For{ $n \in \nodes_t$}
		\State Create an arc $a$ emanating from $n$ for each possible value in $\feasibleSet_t(\nodeState(n))$
		\For{ $a\in \outgoing(n)$ }
		\If{there exists node $n'\in \nodes_{t+1}$ with $\nodeState(n') = \transitionF_t(\nodeState(n), \arcValue_a)$}
		\State Direct arc $a$ to node $n'$, i.e., $\child_a = n'$
		\Else
		\State Create  $n'$ in $\nodes_{t+1}$ with $\nodeState(n') = \transitionF_t(\nodeState(n), \arcValue_a)$ and point arc $a$ to $n'$ 
		\EndIf
		\EndFor
		\EndFor
		\If {$|\nodes_{t+1}| > \maxwidth$} 
		\State Apply \algDDMerge($\nodes_{t+1}$) (relaxed DD) or \algDDDiscard($\nodes_{t+1}$) (restricted DD)
		\EndIf
		\EndFor
		\For{ $n \in \nodes_{\horizon}$}
		\State Create an arc $a$ with  $\parent_a = n$ and $\child_a =\terminalnode$ for each value in $\feasibleSet_\horizon(\nodeState(n))$
		\EndFor
		\State \algReturn\ $\dd$
		\EndProcedure
	\end{algorithmic} 
\end{algorithm}

Note that this procedure can create exact and approximated DDs depending on the maximum width and the procedure to limit the width size. Moreover, this construction mechanism allows the user to obtain a wide range of approximated DDs by simply changing the merge and discarding functions. 

\section{DD-suite: basic features and documentation} \label{sec:software}

DD-suite is an open-source cross-platform code for creating and utilizing DDs to solve discrete optimization problems. Both implementations (C++ and Python) share a similar structure, with identical file names, classes, and function names, albeit with minor differences due to language-specific semantics. We opt for these two languages because they are widely used in the mathematical programming community, and many open-source and commercial solvers support them (e.g., Gurobi \cite{gurobi} and CPLEX \cite{cplex}).

While the two DD-suite implementations are equivalent, we develop each one with consideration for users' potential usage and coding capabilities. We consider the Python version to be a good starting point for people new to DDs and for prototyping purposes. On the other hand, the C++ version could be preferable for researchers interested in computational performance and extending DD-suite capabilities. Nonetheless, users can use any version depending on their preferences and expertise. 

We now detail the basic features of DD-suite, that is, any usage that requires minimal intervention in the existing code. The primary purpose of DD-suite is to enable users to construct and utilize DDs to solve discrete optimization problems easily. Thus, the following sections explain how to specify a problem, build a DD, and obtain optimality bounds. We illustrate these basic features using the Python version for ease of exposition and use \texttt{typewriter} font to refer to code classes and functions. The last section discusses the code documentation and presents a webpage where users can find updated information and step-by-step tutorials.



\subsection{Problem specification}

We now explain how to specify and implement a discrete optimization problem in DD-suite, which requires a recursive formulation similar to the one introduced in Section \ref{sec:background}. Specifically, DD-suite uses \textit{abstract classes} and \textit{inheritance} to specify new problems, allowing us to clearly instruct users on how to implement their recursive formulation and other problem specifications. Thus, the code includes a class called \texttt{AbstractProblem} that serves as the foundation for any problem we want to specify. 

Users have to create a new class that inherits from \texttt{AbstractProblem} and complete the functions that detail their recursive model to introduce a new problem to DD-suite. Similarly to the recursive formulation \eqref{eq:bellman}, the constructor of the new problem class receives the initial state and the variables of the problem. DD-suite accepts any data structure as a state, including pre-defined (e.g., \texttt{int}, \texttt{vector}) and user-specified structures (e.g., \texttt{my\_state}), as long as the user utilizes the same data structure for all states for a particular DD. Additionally, variables can have any bounded integer domain, which can be distinct from one another (e.g., $x_1 \in \{0,1\}$ and $x_2 \in \{0,1,2\}$). The class constructor can also receive any other type of argument that is needed to define the problem (e.g., items' weights in the knapsack problem). 

Next, users need to specify the \texttt{transition\_function()}, which has two primary purposes: (i) detail how to transition from one state to the other and  (ii) identify feasible variable assignments (i.e., implement $\transitionF_t(\cdot)$ and $\feasibleSet_t$, respectively, for any stage $t \in \stages$). Users have complete flexibility when implementing this transition function; however, we recommend a highly efficient implementation, as it is called multiple times during the DD construction procedure and, thus, directly affects the DD construction time.

The constructor and the \texttt{transition\_function()} of the new problem class are the most crucial parts to specify a new problem properly. Our \texttt{AbstractProblem} class also requires users to implement additional functions to ensure that everything works smoothly, which we discuss in more detail in the tutorial provided on our webpage (see Section \ref{sec:webpage}). As an overview, users have to implement the functions \texttt{get\_state\_as\_string()} and \texttt{get\_state\_copy()}, which are used during the DD construction procedure. Additionally, users must specify additional functions when creating relaxed and restricted DDs, which we detail in the following section. 

Note that the environment for problem specification is quite general and can be applied to a wide range of problems. DD-suite provides four problem examples (i.e., knapsack, set cover, maximum independent set, and sequencing problems) for users to experiment with, as well as a SOC knapsack example used in the cutting plane experiments of Section \ref{sec:exp-cuts}. Specifying new problems requires minimal implementation and is quite simple given a recursive formulation, as illustrated in Example \ref{example:software_problem}.

We note that the DD-suite problem specification only describes the feasibility set (i.e., reachable states and feasible variable assignments) and omits the reward function, which is encoded separately (see Section \ref{sec:dd_bounds} for further details). Our implementation intentionally keeps this separation, allowing users to use the same DD graph (i.e., feasibility set) to optimize different objective functions, which is desirable in some DD applications, such as cut generation (see Section \ref{sec:developers}). Nonetheless, users can include reward function information into their problem class if desired (e.g., to guide merging and discarding heuristics) by introducing such information in the class constructor or even adding it to their state definition.

\begin{example}\label{example:software_problem}
Continuing with our running example, Listing \ref{code:knapsack-problem} presents the constructor and transition function implementation for a knapsack problem (Listing \ref{code:knapsackClass} in Appendix \ref{appendix:code} shows a complete implementation of the \texttt{KnapsackProblem} class). The \texttt{KnapsackProblem} class inherits from the \texttt{AbstractProblem} class and implements the recursive model as follows. The constructor receives a single problem-instance object whose attributes bundle the initial state, the variables (i.e., a list of the variable ids and their domain), and the parameters of the knapsack constraint (i.e., items' weights and knapsack capacity). This example uses a single integer as the state representation for the current knapsack load, exactly as in Examples  \ref{example:recursiveKnapsack} and \ref{example:exactDD}.  
The transition function follows the definition of $\transitionF_t(\cdot)$ in Example \ref{example:recursiveKnapsack}: it writes the new state into a scratch buffer and returns whether the resulting load is feasible (i.e., it does not exceed the knapsack capacity).

\begin{lstlisting}[language=Python,
    linerange={4-19},
    caption={Partial knapsack problem class implementation in Python.} ,
    label=code:knapsack-problem]
class KnapsackProblem(AbstractProblem):

    def __init__(self, params: 'KnapsackStructure'):
        super().__init__(params.initial_state, params.variables)
        self.weights: list[int] = params.weights
        self.capacity: int = params.right_side_of_restrictions

    def transition_function(self, previous_state: 'State', variable_index: int, variable_value: int, scratch_state: list) -> bool:

        if variable_value == 0:
            scratch_state[0] = previous_state
            return True

        new_state: int = previous_state + self.weights[variable_index] * variable_value
        scratch_state[0] = new_state
        return new_state <= self.capacity
\end{lstlisting}

To create an instance of the problem, we simply create a class object and initialize the required parameters. Listing \ref{code:knapsack-problem-main} shows the input parameters for our running example with four variables and how to create the problem instance (see Listing \ref{code:knapsackMain} in Appendix \ref{appendix:code} for the complete main code of this running example). \hfill $\square$
\end{example}

\begin{lstlisting}[language=Python,
    linerange={11-27},
    caption={Initialize and create a knapsack problem instance.} ,
    label=code:knapsack-problem-main]
# Maximum width for the approximated DDs
width: int = 2

# Input data for the running example
params = SimpleNamespace(
    initial_state=0,
    variables=[('x_1', [0, 1]), ('x_2', [0, 1]), ('x_3', [0, 1]), ('x_4', [0, 1])],
    weights=[7, 5, 4, 1],
    right_side_of_restrictions=8,
)

# Objective coefficients and sense used later to optimize over the DDs
utility: list[int] = [4, 2, 5, 1]
objective: str = "max"

# Create the knapsack problem used to build the DD
knapsack_problem = KnapsackProblem(params)
\end{lstlisting}

\subsection{DD construction and personalization}

Users can construct a DD by creating an object of the provided \texttt{DD} class, which takes as input an object representing the problem instance. Then, users can call any of the provided functions to create a DD, that is, \texttt{create\_decision\_diagram()} for exact DDs, \texttt{create\_restricted\_decision\_diagram()} for restricted DDs, and \texttt{create\_relax\_priority\_decision\_diagram()} for relaxed DDs. DD-suite includes an alternative relaxed DD constructor, \texttt{create\_relax\_grouping\_decision\_diagram()}, which we further discuss in Section \ref{sec:developer-construction}.

DD-suite utilizes the top-down construction procedure detailed in Algorithm \ref{alg:dd_topdown}. This procedure creates each layer sequentially, following the variable ordering provided by the problem class. Users can control this ordering through the problem class: by default, the variables are processed in the order given to the constructor, but users can define a custom ordering by overriding the \texttt{sort\_variables()} method, which is applied when the problem is constructed with the \texttt{sort} flag enabled. Also, similarly to Algorithm \ref{alg:dd_topdown}, our implementation of the top-down construction procedure is general for all DD types (see Section \ref{sec:developers} for further details) and only differs in whether we keep, discard, or merge nodes when a maximum width is exceeded (for exact, restricted, and relaxed DDs, respectively).

Specifically, function \texttt{create\_decision\_diagram()} constructs an exact DD by simply creating as many nodes as needed in every layer. In contrast, we can create a restricted DD via the command \texttt{create\_restricted\_decision\_diagram()}, which receives a maximum width as input and discards nodes when it is exceeded. Users can specify a discarding heuristic with the \texttt{get\_priority\_for\_discard\_node()} function available in the problem class. This function receives the state of a node and assigns it a priority, where states (i.e., nodes) with higher priority are discarded first.

Similarly, \texttt{create\_relax\_priority\_decision\_diagram()} receives a maximum width as input, and users can specify how to merge nodes in a relaxed DD with functions \texttt{get\_priority\_for\_merge\_nodes()} and \texttt{merge\_operator()} in their problem class. The first function assigns a priority to each node based on its state information. Then, the construction procedure orders all the nodes by priority and repeatedly merges the two highest-priority nodes until the layer fits within the maximum width. The merging operator is specified by the user using the \texttt{merge\_operator()} function, which returns the new merged state. We note that users have complete flexibility in implementing their priority and merging functions, but recall that \texttt{merge\_operator()} must satisfy the properties detailed in Section \ref{sec:relaxedDDs} to ensure that the resulting DD is indeed relaxed. 

While there are many ways to decide which nodes should be discarded or merged, we opt for this priority list type of implementation due to its simplicity and flexibility in accommodating a wide range of heuristics. We note that DD-suite can also be adapted to consider alternative methods for determining which nodes to discard or merge, which is further discussed in Section \ref{sec:developer-construction}. 

Once the DD is created, class \texttt{DD} provides the \texttt{reduce\_decision\_diagram()} function, which applies the reduction procedure mentioned in Section \ref{sec:dd-construction}. Additionally, we can obtain the construction and reduction times using \texttt{get\_building\_time()} and \texttt{get\_reduction\_time()}, respectively, as well as graph information, including the number of nodes and arcs (i.e., \texttt{get\_decision\_diagram().get\_node\_count()} and \texttt{get\_decision\_diagram().get\_arc\_count()}, respectively). Lastly, users can export the DD graph with the \texttt{export\_graph\_file()} command and visualize the resulting DD using, for instance, the yEd software\footnote{Free graph editor software \url{https://www.yworks.com/products/yed}. }. We refer the reader to the DD-suite webpage for further instructions on DD visualization.

\begin{lstlisting}[language=Python,
    linerange={33-42},
    caption={Construct an exact DD, and print relevant information.},
    label=code:create-dd-main]
dd_exact = DD(knapsack_problem)
dd_exact.create_decision_diagram()
dd_exact.reduce_decision_diagram()

# Print DD construction information
print("-- Exact DD information --")
print("\tConstruction time (sec): ", dd_exact.get_building_time())
print("\tReduction time (sec): ", dd_exact.get_reduction_time())
print("\tNumber of nodes: ", dd_exact.get_decision_diagram().get_node_count())
print("\tNumber of arcs: ", dd_exact.get_decision_diagram().get_arc_count())
\end{lstlisting}

\begin{example}
Following up on our running example, Listing \ref{code:create-dd-main} illustrates the command to create and reduce an exact DD for our knapsack instance. Similar commands can be used to create relaxed and restricted DDs, which are available in Appendix \ref{appendix:code}, Listing \ref{code:knapsackMain}. The example code also shows how to obtain the building and reduction times, as well as some key features of the DD graph (i.e., number of nodes and arcs). Lastly, the code illustrates how to export the DD graph into a \texttt{.gml} file using the \texttt{export\_graph\_file()} command, which can then be visualized using the free software yEd. Figure \ref{fig:yEd_DDgraphs} in Appendix \ref{appendix:yEd_graphs} shows the resulting visualizations of such software for all three DD types, which are equivalent to the DDs illustrated in Figures \ref{fig:exactDDKnapsack} and \ref{fig:approxDDKnapsack}.
\end{example}

\subsection{Obtain primal and dual bounds} \label{sec:dd_bounds}

As discussed in Section \ref{sec:background}, we can obtain an optimal solution from a DD using a longest-path algorithm over the graph where arc lengths correspond to the reward function. In the case of approximated DDs, such a procedure provides either primal or dual bounds, depending on whether the DD is restricted or relaxed, respectively. DD-suite includes this optimization feature for cases where the objective function is linear, a common case in applications \cite{castro2022decision,van2024introduction}.

The \texttt{ShortestLongestPath} class provides algorithms for finding the longest or shortest path, depending on whether the problem is maximization or minimization. Its implementation is a simple Bellman--Ford shortest-path algorithm \cite{ahuja1993network} that leverages the structure of the DD graph (i.e., a layered directed acyclic graph) to improve computational efficiency. To use this procedure, users must create an object of the \texttt{ShortestLongestPath} class that receives as input a \texttt{DD} object. Then, users set the arc lengths and the optimization sign (i.e., min or max) with the \texttt{set\_parameters()} function and optimize with the \texttt{solve()} command.  This command returns the solution in a data structure called \texttt{PathStructureSolution}, which includes the optimal value (\texttt{value}), a human-readable description of the optimal path (\texttt{path\_print}), and the list of arcs that form such a path (\texttt{path\_arcs}). The running time of the algorithm can be retrieved separately with the \texttt{get\_time()} command.

We note that each object of \texttt{ShortestLongestPath} is associated with a single DD, but the objective function parameters can be updated as many times as the user wants. This characteristic is particularly useful when using a DD as a sub-routine of a more complex optimization approach, such as Benders decomposition, column generation, or cutting planes (see Section \ref{sec:cutting-planes}).

\begin{example}
    Listing \ref{code:bounds-dd-main} illustrates how to use the \texttt{ShortestLongestPath} class to find an optimal solution over the exact DD for our running example (equivalent code for restricted and relaxed DDs can be found in Appendix \ref{appendix:code}, Listing \ref{code:knapsackMain}). The class object receives the  DD and then updates its cost function with a reward vector (i.e., \texttt{utility}) and the optimization type (i.e., \texttt{objective = "max"}). Lastly, the code executes the longest-path procedure using the \texttt{solve()} command and illustrates how to obtain information about the solution, including its value and the optimal path, as well as the algorithm's execution time via \texttt{get\_time()}.
\end{example}

\begin{lstlisting}[language=Python,
    linerange={81-88},
    caption={Obtain optimal solution and print relevant information.},
    label=code:bounds-dd-main]
# Optimal solution over the exact DD (longest path)
longest_path = ShortestLongestPath(dd_exact)
longest_path.set_parameters(utility, objective)
answer = longest_path.solve()

print("The optimal value is: ", answer.value)
print("Longest path: ", answer.path_print)
print("Solution time (sec): ", longest_path.get_time())
\end{lstlisting}

\subsection{Documentation, examples, and webpage} \label{sec:webpage}

To facilitate DD-suite usage, our code provides thorough documentation of all classes and functions. This documentation follows standard software engineering documentation practices (i.e., descriptive docstrings for every class, method, and function, including their parameters and return values) and is included in the foundation code and all the provided examples. Moreover, DD-suite follows the Clean Code principles \cite{martin2009clean} (e.g., meaningful names, small functions and classes, etc.), so it is user-friendly and self-explanatory. We further detail the advantages of following such software engineering principles in Section \ref{sec:developers}. 

Users can run the examples described above using the provided data instances and one of the main files provided to reproduce our experiments (see Section \ref{sec:experiments} for further details) or create their own main files to try out DD-suite functionalities. Moreover, these examples include automatic testing, so users can easily check if the code still runs smoothly if they make changes to it (see Section \ref{sec:developers} for further details).

In addition, we have set up a webpage for DD-suite\footnote{\url{https://dd-suite.dev/}}. The site contains information on how to set up the Python and C++ versions and run the code to reproduce our experiments. Moreover, we provide a step-by-step tutorial that focuses on the basic features detailed in this section and some advanced features (i.e., cutting planes, see Section \ref{sec:cutting-planes}) for both the Python and C++ implementations. The webpage also includes essential information for developers, such as how to run and create new tests and a description of the most important classes of the software. 

\section{DD-suite for developers: advanced usage and extensions} \label{sec:developers}

As previously mentioned, DD-suite allows users to specify discrete optimization problems easily, create DDs for them, and obtain optimality bounds with just a few lines of code. However, this is just the tip of the iceberg given the diverse uses of DDs in the literature \cite{castro2022decision,van2024introduction}. This section presents the main software characteristics of DD-suite that make it an ideal code base for DD-based algorithms and applications. We then discuss how users can extend and apply DD-suite for their research projects, including our DD-based cutting plane implementation as an illustrative example.

\subsection{Software characteristics and tests} \label{sec:dev-characteristics}

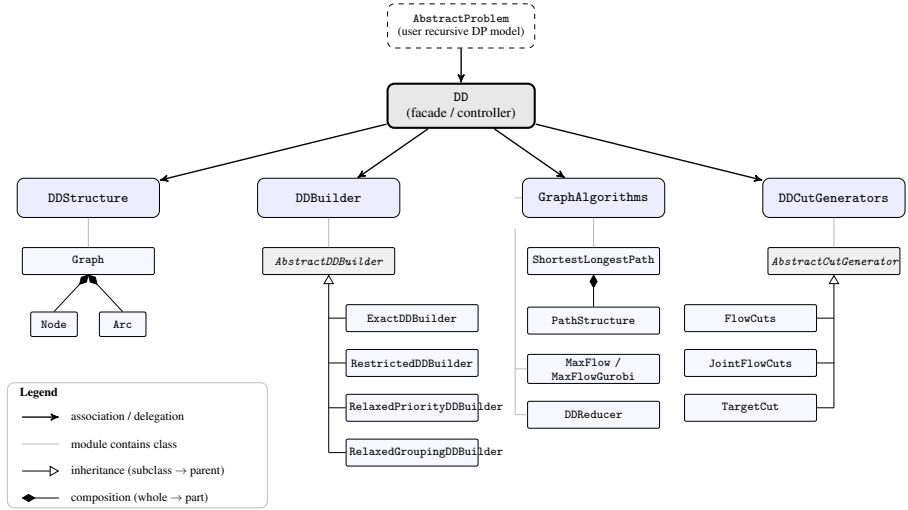
\begin{figure}[t]
    \centering
    \resizebox{\textwidth}{!}{%
    \begin{tikzpicture}[
        >=stealth',
        module/.style={draw, rounded corners, align=center, font=\footnotesize\ttfamily, minimum height=0.8cm, minimum width=2.9cm, fill=blue!7},
        facade/.style={draw, rounded corners, align=center, font=\footnotesize, minimum height=0.9cm, minimum width=3cm, fill=gray!18, very thick},
        io/.style={draw, dashed, rounded corners, align=center, font=\scriptsize, minimum height=0.9cm, minimum width=3cm},
        cls/.style={draw, rounded corners=1pt, align=center, font=\scriptsize\ttfamily, text width=2.55cm, minimum height=0.55cm, fill=blue!3, inner sep=2pt},
        smallcls/.style={draw, rounded corners=1pt, align=center, font=\scriptsize\ttfamily, minimum width=0.95cm, minimum height=0.5cm, fill=blue!3, inner sep=1pt},
        abs/.style={cls, font=\scriptsize\ttfamily\itshape, fill=gray!12},
        assoc/.style={->, thick},
        member/.style={-, gray!55},
        inh/.style={-open triangle 60},
        comp/.style={diamond-},
    ]
        \node[io] (problem) {{\ttfamily AbstractProblem}\\ (user recursive DP model)};
        \node[facade, below=0.7cm of problem] (dd) {{\ttfamily DD}\\ (facade / controller)};
        \draw[assoc] (problem) -- (dd);

        \node[module, below=1.0cm of dd, xshift=-7.6cm]  (mStruct) {DDStructure};
        \node[module, below=1.0cm of dd, xshift=-2.7cm]  (mBuild)  {DDBuilder};
        \node[module, below=1.0cm of dd, xshift=2.7cm]   (mAlg)    {GraphAlgorithms};
        \node[module, below=1.0cm of dd, xshift=7.6cm]   (mCut)    {DDCutGenerators};
        \draw[assoc] (dd) -- (mStruct);
        \draw[assoc] (dd) -- (mBuild);
        \draw[assoc] (dd) -- (mAlg);
        \draw[assoc] (dd) -- (mCut);

        \node[cls, below=0.6cm of mStruct] (graph) {Graph};
        \node[smallcls, below=0.75cm of graph, xshift=-0.7cm] (node) {Node};
        \node[smallcls, below=0.75cm of graph, xshift=0.7cm]  (arc)  {Arc};
        \draw[member] (mStruct) -- (graph);
        \draw[comp] (graph.south) -- (node.north);
        \draw[comp] (graph.south) -- (arc.north);

        \node[abs, below=0.6cm of mBuild] (absB) {AbstractDDBuilder};
        \node[cls, below=0.6cm of absB, xshift=1.7cm] (exact) {ExactDDBuilder};
        \node[cls, below=0.35cm of exact] (restr) {RestrictedDDBuilder};
        \node[cls, below=0.35cm of restr] (relP)  {RelaxedPriorityDDBuilder};
        \node[cls, below=0.35cm of relP]  (relG)  {RelaxedGroupingDDBuilder};
        \draw[member] (mBuild) -- (absB);
        \coordinate (bbBot) at (absB.south |- relG.west);
        \draw[inh] (bbBot) -- (absB.south);              
        \draw (exact.west) -- (exact.west -| absB.south); 
        \draw (restr.west) -- (restr.west -| absB.south);
        \draw (relP.west)  -- (relP.west  -| absB.south);
        \draw (relG.west)  -- (relG.west  -| absB.south);

        \node[cls, below=0.6cm of mAlg] (slp) {ShortestLongestPath};
        \node[cls, below=0.65cm of slp] (paths) {PathStructure};
        \node[cls, below=0.4cm of paths] (maxflow) {MaxFlow /\\ MaxFlowGurobi};
        \node[cls, below=0.35cm of maxflow] (reducer) {DDReducer};
        \draw[member] (mAlg) -- (slp);
        \draw[comp] (slp.south) -- (paths.north);
        \coordinate (baT) at ([xshift=-1.6cm,yshift=-0.22cm]mAlg.south);
        \draw[member] (baT) -- (baT |- reducer.west);
        \draw[member] (mAlg.west)    -- (mAlg.west    -| baT);
        \draw[member] (maxflow.west) -- (maxflow.west -| baT);
        \draw[member] (reducer.west) -- (reducer.west -| baT);

        \node[abs, below=0.6cm of mCut] (absC) {AbstractCutGenerator};
        \node[cls, below=0.6cm of absC, xshift=-1.7cm] (flow)  {FlowCuts};
        \node[cls, below=0.35cm of flow]  (joint) {JointFlowCuts};
        \node[cls, below=0.35cm of joint] (target){TargetCut};
        \draw[member] (mCut) -- (absC);
        \coordinate (bcBot) at (absC.south |- target.east);
        \draw[inh] (bcBot) -- (absC.south);
        \draw (flow.east)   -- (flow.east   -| absC.south);
        \draw (joint.east)  -- (joint.east  -| absC.south);
        \draw (target.east) -- (target.east -| absC.south);

        \coordinate (legO) at ([xshift=-1.35cm,yshift=-2.7cm]graph.south);
        \begin{scope}[shift={(legO)}]
            \draw[rounded corners, draw=gray!50] (-0.3,0.5) rectangle (5.0,-2.05);
            \node[anchor=west, font=\scriptsize\bfseries] at (-0.15,0.28) {Legend};
            \draw[assoc]  (0,-0.2) -- (0.75,-0.2);   \node[anchor=west, font=\scriptsize] at (0.9,-0.2)  {association / delegation};
            \draw[member] (0,-0.75) -- (0.75,-0.75); \node[anchor=west, font=\scriptsize] at (0.9,-0.75) {module contains class};
            \draw[inh]    (0,-1.3) -- (0.75,-1.3);   \node[anchor=west, font=\scriptsize] at (0.9,-1.3)  {inheritance (subclass $\rightarrow$ parent)};
            \draw[comp]   (0,-1.85) -- (0.75,-1.85); \node[anchor=west, font=\scriptsize] at (0.9,-1.85) {composition (whole $\rightarrow$ part)};
        \end{scope}
    \end{tikzpicture}%
    }
    \caption{Class-level architecture of DD-suite.}
    \label{fig:architecture}
\end{figure}

DD-suite is designed so that researchers can add new procedures, such as algorithms running on top of a DD or new DD construction mechanisms, without rewriting the existing code. We achieve this by following Clean Code software engineering principles \cite{martin2009clean}: meaningful names, small single-purpose functions, and, most importantly, a modular architecture where every procedure is encapsulated in its own component. The code is organized into independent modules, namely the DD data structures (\texttt{DDStructure}), the construction procedures (\texttt{DDBuilder}), the algorithms that interact directly with the graphical structure (e.g., shortest/longest path and reduction), and advanced algorithms for optimization purposes (e.g., cutting planes). The \texttt{DD} class acts as a controller (i.e., a facade) that coordinates these modules and exposes a clean interface to the user. As a consequence, the components are decoupled: changing how a DD is constructed does not affect the algorithms that run on top of it, and vice versa, as long as the shared \texttt{Graph} interface exposed by \texttt{DDStructure} is respected. This separation is what makes DD-suite easy to extend, as we illustrate in the following sections. Figure \ref{fig:architecture} summarizes these modules and how they interact through the \texttt{DD} facade, together with the main classes inside each module and their relationships. In particular, the figure shows how the data structures are composed, where a \texttt{Graph} is built from \texttt{Node} and \texttt{Arc} objects and a \texttt{ShortestLongestPath} from \texttt{PathStructure} objects, as well as the abstract classes that new procedures specialize, which we discuss in the following sections.

A key ingredient that enables modularity is that all construction procedures share a single top-down construction template, implemented in the \texttt{AbstractDDBuilder} class following Algorithm \ref{alg:dd_topdown}. This template builds the diagram layer by layer and exposes two extension points (i.e., abstract methods) that each DD type implements: one invoked after every layer is completed and one invoked at the end of the construction. The exact, restricted, and relaxed builders are all concrete subclasses of \texttt{AbstractDDBuilder} that differ only in these two methods (i.e., whether they keep, discard, or merge nodes when the maximum width $\maxwidth$ is exceeded), while reusing the rest of the construction machinery unchanged.

To safeguard this code base against future modifications, DD-suite includes an extensive automatic test suite for both implementations, using the \texttt{unittest}/\texttt{pytest} frameworks in Python and Google Test in C++. At the time of writing, the suite comprises more than $230$ tests in Python and $260$ in C++, and this number keeps growing as new problems and procedures are added. The tests cover each example problem and DD type (exact, restricted, and relaxed, with and without custom variable ordering), the graph operations (e.g., reduction and graph equivalence), and the cut algorithms. This testing environment serves a double purpose: it protects the existing functionality, and it gives developers a safety net to experiment with new ideas with the guarantee that they are not silently breaking the original code. We provide instructions on how to run and write new tests on the DD-suite webpage.

\subsection{Extending the DD construction mechanism} \label{sec:developer-construction}

As previously mentioned, the top-down construction procedure is implemented in the \texttt{AbstractDDBuilder} class, which builds the diagram following Algorithm \ref{alg:dd_topdown} and delegates the type-specific behavior to two methods: \texttt{\_specific\_end\_of\_layer \_function()}, called after each layer is built, and \texttt{\_specific\_end\_of\_construction \_function()}, called once after the last layer. To incorporate a new construction mechanism, a developer only needs to create a new subclass of \texttt{AbstractDDBuilder} and implement these two methods, reusing all the underlying machinery.

We illustrate this extension mechanism with the alternative grouping-based relaxation \texttt{create\_relax\_grouping\_decision\_diagram()} routine. 
Recall that the default relaxed DD builder (i.e., \texttt{RelaxedPriorityDDBuilder}) merges the two highest-priority nodes. In contrast, the grouping variant (i.e., \texttt{RelaxedGroupingDDBuilder}) merges several nodes at once, as described in Algorithm \ref{alg:dd_grouping_merge}: it sorts the nodes of the layer by merge priority and, in successive passes, merges every group of consecutive nodes whose priority gap lies within a tolerance $\delta$. The tolerance starts at $\delta=0$, so that only nodes with equal priority are merged, and grows after each pass until the layer fits within the maximum width $\maxwidth$; to avoid passes that would merge nothing, $\delta$ jumps to the smallest priority gap still present in the layer whenever that gap is larger. Note that this procedure is a drop-in replacement for the \algDDMerge\ routine of Algorithm \ref{alg:dd_topdown} (line 13) and reuses the same user-defined priority and merge operators. 

\begin{algorithm}[b]
	\caption{Grouping-based Merge of a Node Layer} \label{alg:dd_grouping_merge}
	\begin{algorithmic}[1]
		\Procedure{\algDDMergeGroup}{$\nodes_t$, $\maxwidth$}
		\State Sort $\nodes_t$ in non-decreasing order of the merge priority $\priorityF(\cdot)$
		\State $\delta \gets 0$ \Comment{tolerance: merge only nodes with equal priority}
		\While{$|\nodes_t| > \maxwidth$}
		\State Partition $\nodes_t$ into maximal groups of consecutive nodes whose priority gap is at most $\delta$
		\For{each group $\mathcal{G}$ with $|\mathcal{G}| > 1$ and while $|\nodes_t| > \maxwidth$}
		\State Create node $\tilde{n}$ with $\nodeState(\tilde{n}) = \bigmerge_{n \in \mathcal{G}} \nodeState(n)$
		\State Redirect to $\tilde{n}$ every arc pointing to a node of $\mathcal{G}$, and replace $\mathcal{G}$ by $\tilde{n}$ in $\nodes_t$
		\EndFor
		\State $\delta \gets \max\{ \delta + 1, \; \text{smallest priority gap in } \nodes_t \}$ \Comment{skip passes that merge nothing}
		\EndWhile
		\State \algReturn\ $\nodes_t$
		\EndProcedure
	\end{algorithmic}
\end{algorithm}

Listing \ref{code:dd-builder-extension} shows the skeleton of this new builder: implementing the alternative relaxation required overriding a single layer-level step (the grouped merge), while the rest of the construction (i.e., the node creation, the priority and merge operators defined in the problem class, and the final pruning) is inherited unchanged. The same pattern can be used to prototype other construction strategies discussed in the literature, such as iterative refinement or separation-based constructions \cite{cire2014separation,van2022graph}.

\begin{lstlisting}[language=Python,
    caption={Skeleton of the grouping-relaxation builder (abbreviated).},
    label=code:dd-builder-extension]
class RelaxedGroupingDDBuilder(AbstractDDBuilder):

    def __init__(self, problem, max_width):
        super().__init__(problem)
        self._max_width = max_width

    def _specific_end_of_layer_function(self):
        # Called after every layer: shrink it if it exceeds the maximum width.
        if len(self.graph.structure[-1]) > self._max_width:
            self._merge_nodes_when_width_is_greater_than_w()
            self._delete_nodes_current_layer(self.graph.actual_layer)

    def _specific_end_of_construction_function(self):
        # Called once after the last layer: prune dead-ends and renumber nodes.
        self._bottom_up_pruner()
        self.adjust_node_number()
\end{lstlisting}
    
The key difference with the default relaxed builder lies in the end-of-layer step: the grouped merge collapses many nodes in a single pass, which is cheaper on wide layers, at the cost of coarser control over which nodes are merged. All the remaining machinery (node and arc creation, the transition function, the user-defined priority and merge operators, and the final pruning) is reused from \texttt{AbstractDDBuilder} without modification.

\subsection{Implementing new algorithms} \label{sec:dev-algorithms}

Beyond constructing DDs, many DD-based techniques run algorithms on top of the resulting diagram (e.g., shortest paths for bounds, max-flow for separation, or reduction). DD-suite supports this usage directly: once a DD is built, the \texttt{DD} class exposes the underlying graph through \texttt{get\_decision\_diagram()} (and a deep copy through \texttt{get\_decision\_diagram\_copy()}), so that any graph algorithm can be implemented against the same \texttt{Graph}, \texttt{Node}, and \texttt{Arc} data structures. The algorithms included in DD-suite are concrete examples of this pattern and live in the \texttt{GraphAlgorithms} module: the \texttt{ShortestLongestPath} class that computes optimality bounds (Section \ref{sec:dd_bounds}), the \texttt{MaxFlow} and \texttt{MaxFlowGurobi} routines used by the cut generators, and the \texttt{DDReducer} that performs the reduction procedure.

A design decision that is particularly convenient for developing new algorithms is the separation between the feasibility set and the reward function (see Section \ref{sec:dd_bounds}). Because the DD encodes only the reachable states and feasible assignments, the same graph can be reused to optimize over many different objective functions without rebuilding it. This is exactly what advanced procedures (e.g., Benders decomposition, column generation, and cutting planes) require: a DD is queried repeatedly with changing objectives.

\subsection{Example: A cutting plane implementation} \label{sec:cutting-planes}

To demonstrate how DD-suite serves as a building block for more sophisticated algorithms, we implement a DD-based cutting plane procedure and embed it into a state-of-the-art mixed-integer programming (MIP) solver. This functionality is encapsulated in the \texttt{DDCutGenerators} module, which follows the same abstract-class design used throughout the code. The \texttt{AbstractCutGenerator} class, shown in Listing \ref{code:cut-generator}, defines the interface that every cut generator must implement: given a (fractional) point, \texttt{generate\_cut()} attempts to separate it and returns whether a violated inequality was found, while \texttt{get\_cut()} returns the coefficients and right-hand side of the resulting inequality. This uniform interface lets the host solver treat all cut families interchangeably.

\begin{lstlisting}[language=Python,
    caption={The \texttt{AbstractCutGenerator} interface (abbreviated).},
    label=code:cut-generator]
class AbstractCutGenerator(ABC):

    def __init__(self, tolerance=1e-4):
        self.tolerance = tolerance

    @abstractmethod
    def generate_cut(self, x_values, verbose=False) -> bool:
        # Separate the point x_values; return True if a violated cut was found.
        ...

    @abstractmethod
    def get_cut(self) -> tuple[list[float], float]:
        # Coefficients and right-hand side of the last separated inequality.
        ...

    @abstractmethod
    def get_name(self) -> str: ...

    @abstractmethod
    def get_time(self) -> float: ...
\end{lstlisting}

DD-suite provides three families of cut generators built on top of a DD, each as a subclass of \texttt{AbstractCutGenerator}: combinatorial flow cuts (\texttt{FlowCuts}) and their dual max-flow counterpart (\texttt{JointFlowCuts}), both based on the combinatorial cut-and-lift procedure of \cite{castro2022combinatorial}, and target cuts (\texttt{TargetCut}) following \cite{tjandraatmadja2019target}. In addition, the \texttt{CutStrengthening} class strengthens a generated inequality following a disjunctive coefficient strengthening strategy that leverages the DD structure \cite{castro2022combinatorial}. 

Finally, we integrate these generators into the Gurobi \cite{gurobi} B\&B search through a user-cut callback. The \texttt{AbstractProblemGurobiClass} sets up the MIP model and, when DD cuts are enabled, registers a callback that, at the root node, retrieves the current fractional solution (\texttt{cbGetNodeRel}), invokes the active cut generators over the corresponding DDs, and injects any violated inequality back into the solver (\texttt{cbCut}). We emphasize that the entire implementation is built on top of the abstractions described above, with no modification to the DD construction or graph algorithms. 

\section{Numerical Experiments} \label{sec:experiments}

This section presents an empirical evaluation of DD-suite with three main objectives. First, we assess the \emph{computational efficiency} of our cross-platform implementation by comparing the construction time of its Python and C++ versions on a large set of instances and benchmarking them against \texttt{ddo} \cite{gillard2020ddo} (Section \ref{sec:exp-efficiency}). Second, we showcase the \emph{extensibility} of the code by introducing an alternative relaxed-DD construction mechanism and comparing it against the default one (Section \ref{sec:exp-grouping}). Third, we illustrate the \emph{research value} of DD-suite as a building block for DD-based algorithms by implementing a DD-based cutting plane procedure and embedding it into a state-of-the-art MIP solver (Section \ref{sec:exp-cuts}). All scripts needed to reproduce the experiments, including data-visualization routines, are provided with our code\footnote{GitHub repository: \url{https://github.com/MargaritaCastro/dd-suite}}.

\subsection{Experimental setup} \label{sec:exp-setup}

We now describe the experimental setup common to all three experiments, namely the problems and instances, and the computational environment used throughout.

\paragraph{Problems and instances.} We evaluate DD-suite on the five discrete optimization problems included as examples in the code: knapsack, maximum independent set, set cover, sequencing, and SOC knapsack problems. The latter is used only in our cutting plane experiment of Section \ref{sec:exp-cuts}. For each problem, we consider up to two families of instances. The \emph{standard} family collects benchmarks taken from the literature: the knapsack instances of \cite{pisinger2005hard}, the DIMACS clique instances for the independent set problem \cite{johnson1996dimacs}, the SOC knapsack benchmark (SOC-K) of \cite{atamturk2010conic} composed of $90$ instances, and, for the sequencing problem, $40$ instances whose setup time matrices are taken from the ATSP set of TSPLIB \cite{reinelt1991tsplib}. The \emph{custom} family consists of randomly generated instances created using the generators in \texttt{DataInstances/}, which allow us to control problem size and difficulty (e.g., number of variables, density, and seed). All instances of the set cover problem are randomly generated. In total, the benchmark comprises $608$ instances across all problems and families. Appendix \ref{appendix:data-instances} details each problem, the specific instance sets used, and how they were generated.

\paragraph{Computational environment.} All experiments were run on a high-performance computing cluster managed by SLURM, on compute nodes equipped with two Intel\textsuperscript{\textregistered} Xeon\textsuperscript{\textregistered} E5-2630 v4 processors ($2.20$ GHz, $20$ physical cores) and $755$ GB of RAM, running Rocky Linux $9.7$. Each instance was solved on a single thread to enable a fair comparison across implementations. Unless stated otherwise, we build approximated DDs with a maximum width $\maxwidth \in \{500, 1000, 2000, 5000, 10000, 20000\}$ and impose a time limit of five minutes per DD construction; runs exceeding this limit are reported as timeouts. The cutting plane experiments in Section \ref{sec:exp-cuts} use a one-hour time limit, a common choice in the MIP literature. We use Gurobi $13$ \cite{gurobi} both as a baseline MIP solver and as the host solver for our cutting plane procedure.

\subsection{Cross-language efficiency of DD-suite} \label{sec:exp-efficiency}

Our first experiment evaluates the computational performance of the two DD-suite implementations and positions them relative to an existing high-performance tool. For every instance and DD type (exact, restricted, and relaxed), we construct the corresponding diagram in Python and C++ and measure its construction time. The two DD-suite implementations share the same algorithms and produce identical diagrams and optimality bounds; thus, comparing them isolates the overhead of the programming language rather than algorithmic differences. We verified this equivalence by checking that the optimization value, the number of nodes, and the number of arcs coincide across both languages on every configuration, for all DD types, and that the value obtained from the exact DDs matches the optimal value reported by Gurobi (see Appendix \ref{appendix:exact-values}). We further benchmark DD-suite against \texttt{ddo} \cite{gillard2020ddo}, a state-of-the-art compiled (Rust) framework for DD-based optimization, to gauge how our implementation compares against a specialized, performance-oriented tool.

\begin{figure}[b]
    \centering
    \includegraphics[width=0.46\textwidth]{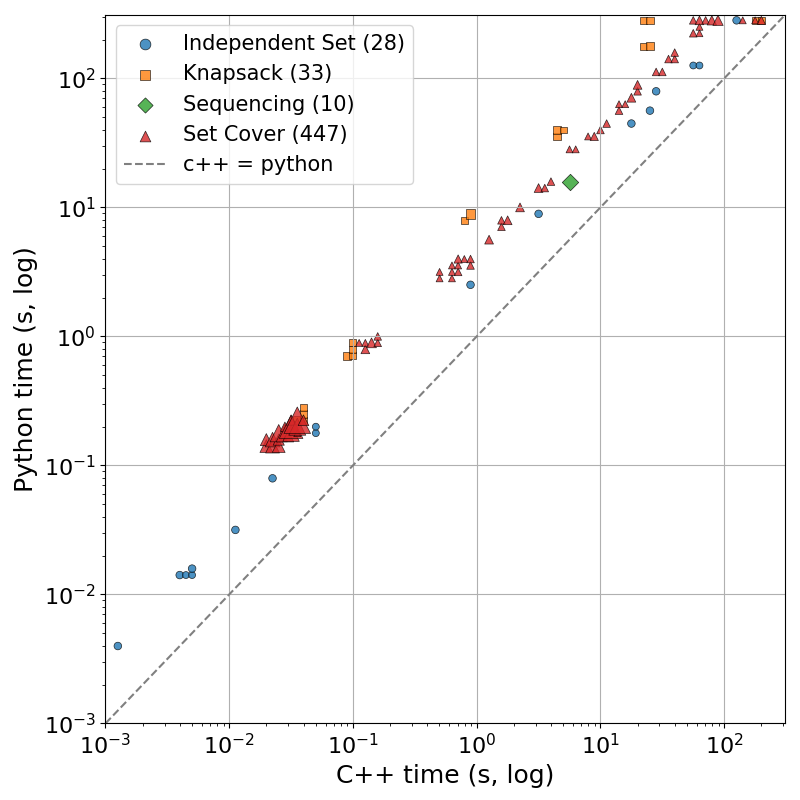}
    \includegraphics[width=0.46\textwidth]{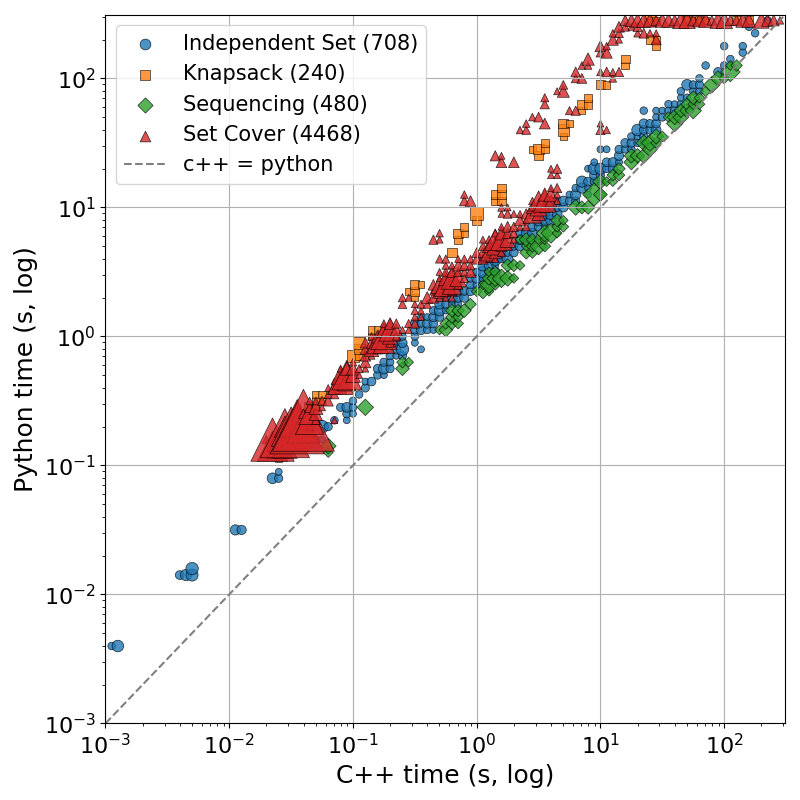}\\[1ex]
    \includegraphics[width=0.46\textwidth]{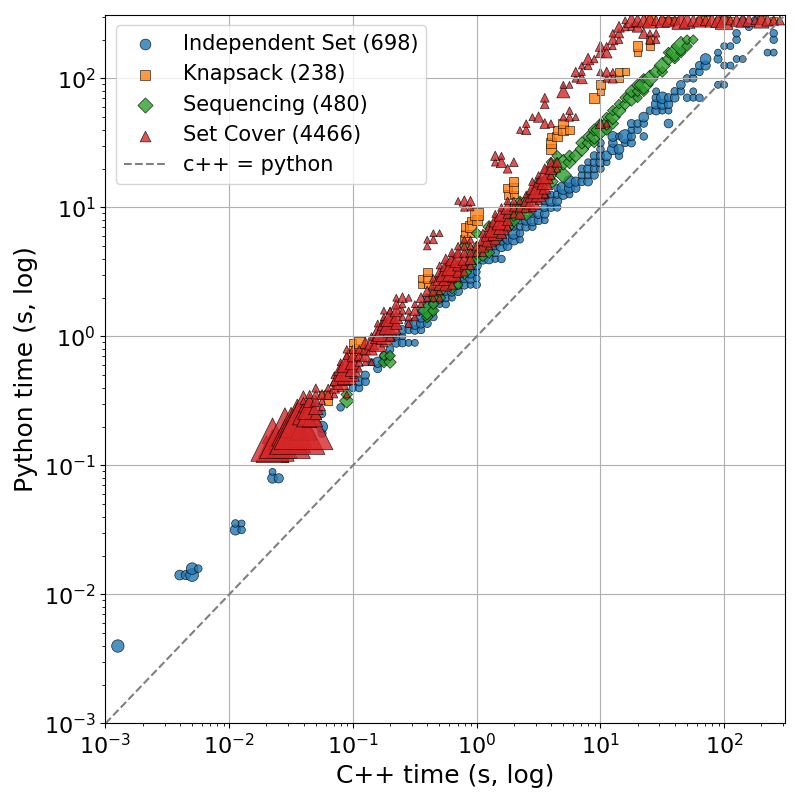}
    \caption{C++ versus Python construction time, per DD type.}
    \label{fig:construction-time}
\end{figure}

Figure \ref{fig:construction-time} compares the construction time of C++ against Python on a logarithmic scale, with one panel per DD type (exact, restricted, and relaxed) and colors denoting the problem type. To avoid overlaps, instances are binned by their (rounded) construction times, so each marker aggregates all instances that fall on approximately the same coordinate; the size of a marker is proportional to the number of instances it represents (i.e., a larger marker indicates that more instances accumulate at that position). Since the $x$-axis corresponds to C++ and the $y$-axis to Python, points above the diagonal indicate instances for which C++ is faster. The C++ implementation is consistently faster than the Python implementation across all problems and DD types, with an average speedup of $5$ to $6$ times. Table \ref{tab:efficiency} summarizes these results, reporting the average construction-time ratio for each DD type, aggregated over all problems and widths, where a construction-time ratio of $k$ means that the slower implementation takes $k$ times longer. 

\begin{table}[tb]
\centering
\caption{Average ratio of construction times per DD type.}
\label{tab:efficiency}
\begin{tabular}{llc}
\toprule
Comparison & DD type & Average construction-time ratio \\
\midrule
Python\,/\,C++ & Exact       & $5.9$ \\
Python\,/\,C++ & Restricted  & $5.2$ \\
Python\,/\,C++ & Relaxed      & $5.7$ \\
\midrule
C++\,/\,\texttt{ddo} & Exact         & $2.1$ \\
C++\,/\,\texttt{ddo} & Restricted    & $2.3$ \\
C++\,/\,\texttt{ddo} & Relaxed        & $1.8$ \\
\bottomrule
\end{tabular}
\end{table}

\begin{figure}[tb]
    \centering
    \includegraphics[width=0.46\textwidth]{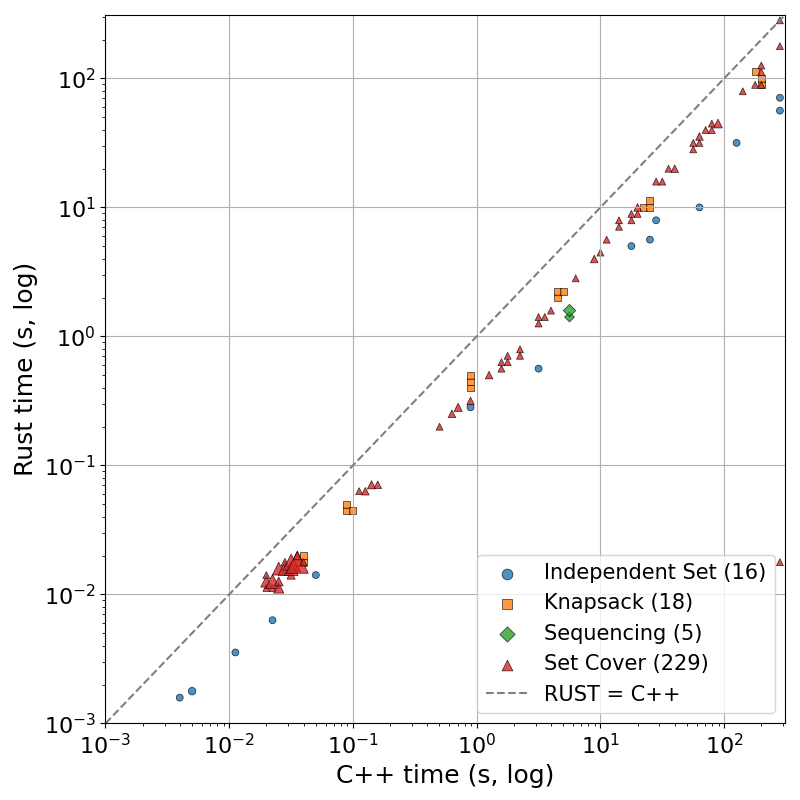}
    \includegraphics[width=0.46\textwidth]{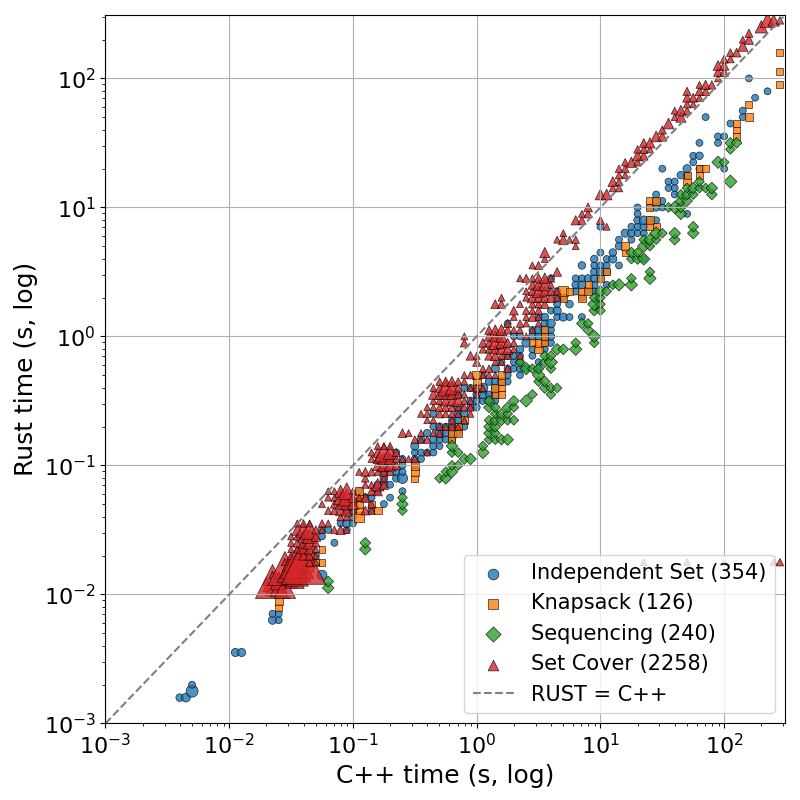}\\[1ex]
    \includegraphics[width=0.46\textwidth]{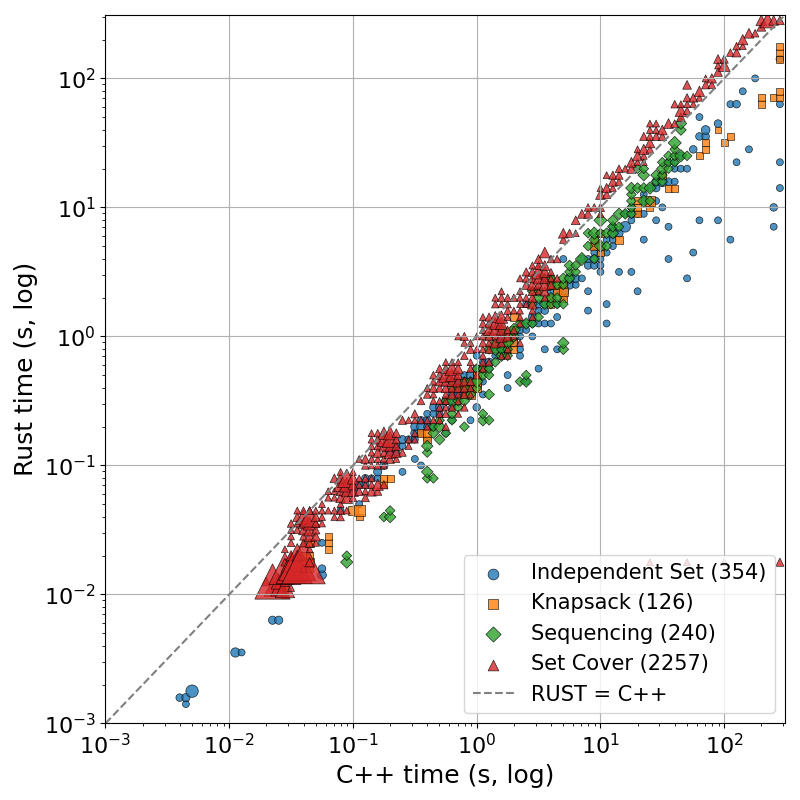}
    \caption{C++ versus \texttt{ddo} construction time, per DD type.}
    \label{fig:construction-time-rust}
\end{figure}

Figure \ref{fig:construction-time-rust} compares DD-suite's C++ implementation against \texttt{ddo}, with one panel per DD type and points below the diagonal indicating that \texttt{ddo} is faster, following the same style as in Figure \ref{fig:construction-time}. We can see that \texttt{ddo} is faster than our C++ implementation with a factor of approximately $2$ for all DD types on average (see Table \ref{tab:efficiency}).
This gap is expected, as \texttt{ddo} is a specialized framework written in Rust and heavily optimized for raw performance, whereas DD-suite prioritizes a clean, extensible, and cross-platform code base. Specifically, DD-suite includes additional data structures for manipulating the underlying graph (e.g., both incoming and outgoing arcs), which makes it amenable to algorithm development but compromises its speed. In contrast, \texttt{ddo} solely focuses on B\&B optimization, so it only includes the necessary structures to favor speed (e.g., only outgoing arcs).


\subsection{Extending DD-suite: an alternative relaxation mechanism} \label{sec:exp-grouping}

This experiment illustrates how DD-suite can be extended within the code with minimal effort. As described in Section \ref{sec:software}, the default relaxed-DD construction merges node pairs following a user-defined priority. We complement it with an alternative routine that merges nodes in groups, which only required overriding a single construction step, as explained in Section \ref{sec:developer-construction}.

\begin{figure}[ht]
    \centering
    \begin{subfigure}[t]{0.49\textwidth}
        \includegraphics[width=1\textwidth]{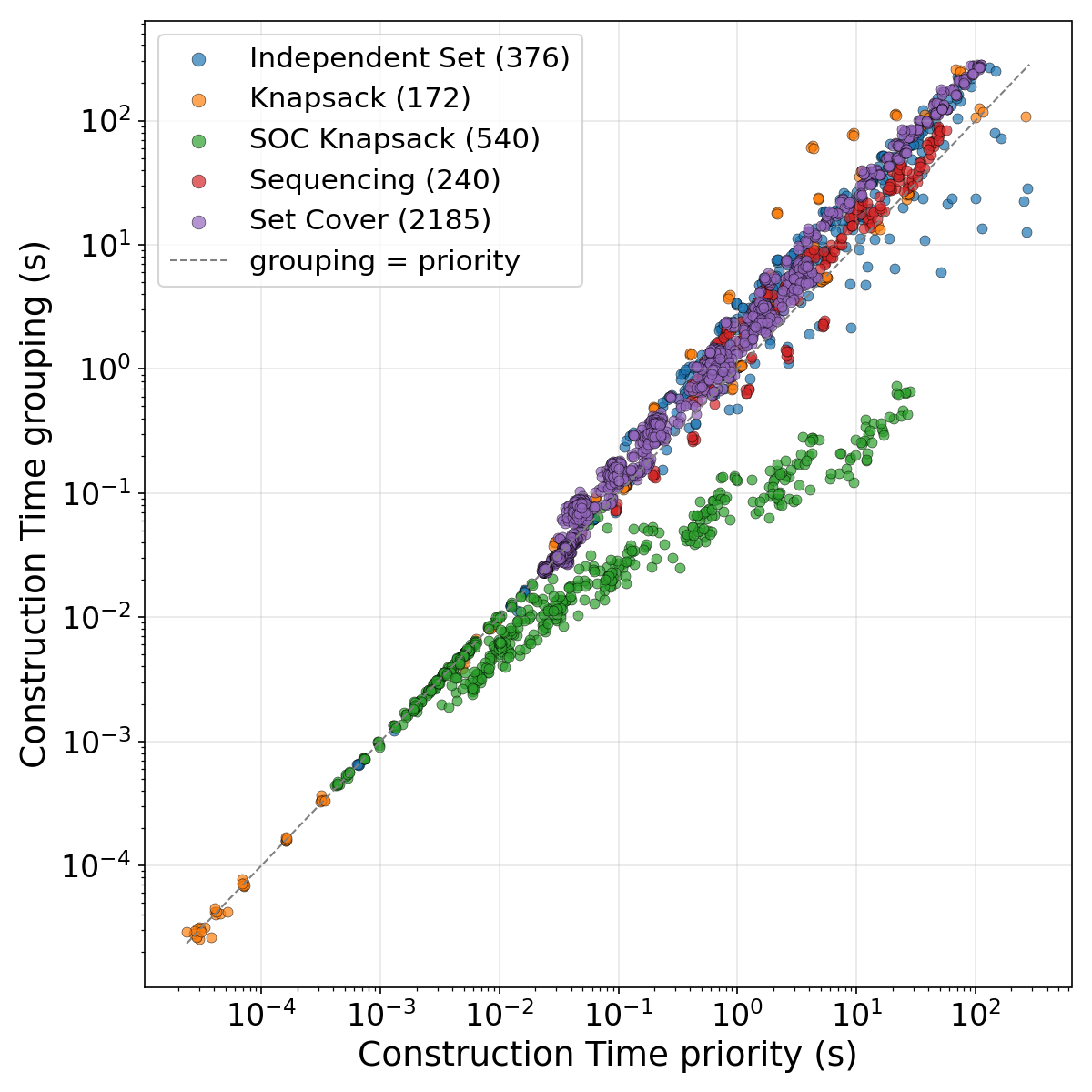}
        \caption{Construction time comparison.}
        \label{fig:grouping-time}
    \end{subfigure}
    ~
    \begin{subfigure}[t]{0.49\textwidth}
        \includegraphics[width=1\textwidth]{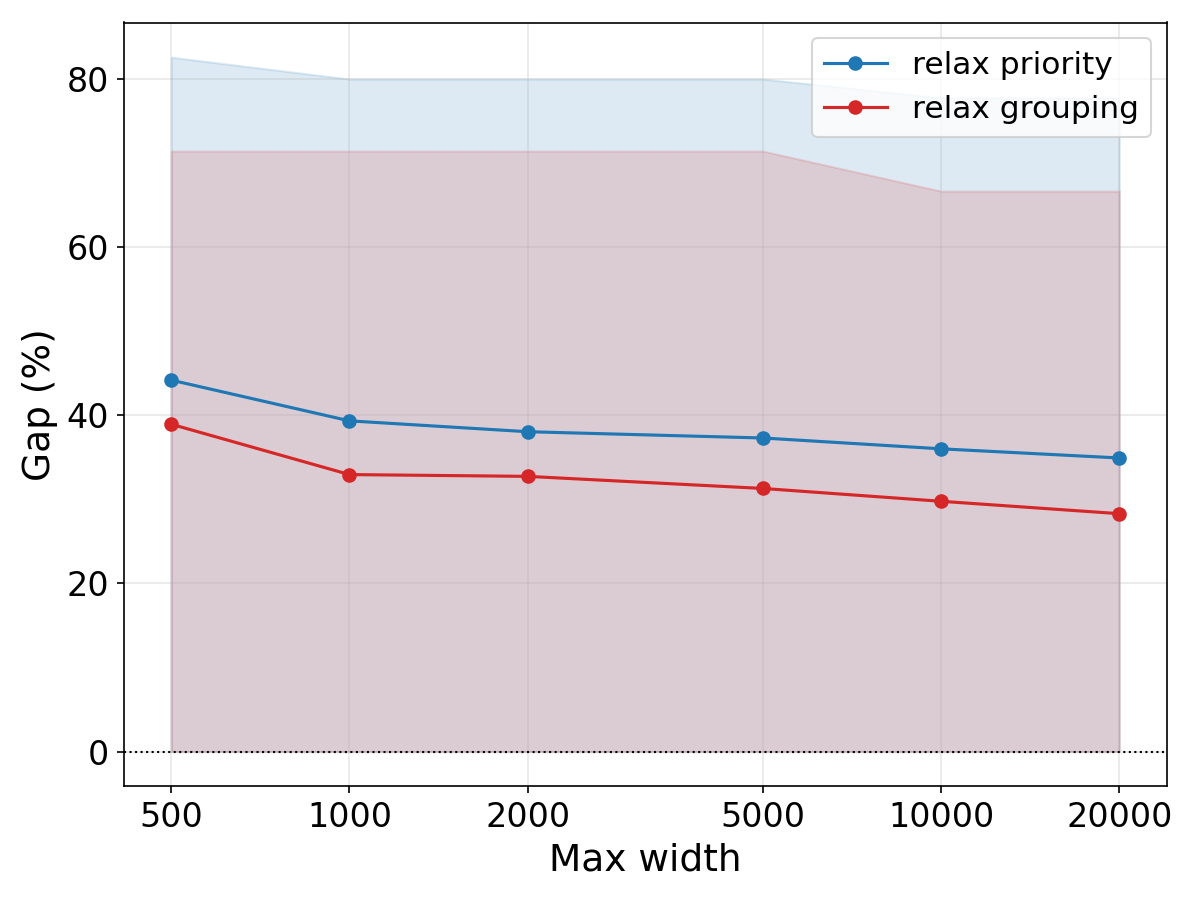}
        \caption{Gap comparison for set cover instances.}
        \label{fig:grouping-gap}
    \end{subfigure}
    \caption{Time and gap comparison between priority and grouping relaxed DDs.}
\end{figure}

Figure \ref{fig:grouping-time} compares the two relaxation mechanisms in terms of construction time, with the priority mechanism on the $x$-axis and the grouping mechanism on the $y$-axis (points below the diagonal indicate that grouping is faster). The grouping variant merges several nodes at once and is therefore especially effective on wide layers. While the two mechanisms have comparable construction times on most instances, grouping is faster on average (with an average construction-time ratio of about $1.8$), with the largest gains concentrated on the instances and widths that produce the widest layers.

Figure \ref{fig:grouping-gap} compares the two mechanisms in terms of bound quality, measured as the optimality gap of the relaxed bound relative to the exact optimum. For each maximum width, the figure reports the average gap (after removing outliers using the IQR rule) and the interquartile range (shaded band) across the set-cover instances for both mechanisms. The two mechanisms yield identical bounds for roughly $60\%$ of the (instance, width) pairs and are comparable on the remainder, with neither dominating the other. The analogous gap plots for all problems and for both the C++ and Python implementations are collected in Appendix \ref{appendix:grouping-gaps}.

We stress that these two construction heuristics were not heavily engineered: our goal is to illustrate that alternative construction mechanisms can be prototyped within DD-suite with minimal effort and that they lead to observable differences in construction time and bound quality, rather than to identify the best possible relaxation. A thorough study of which construction strategy is preferable, and under which conditions, is an interesting future work direction. Overall, the grouping mechanism offers a useful speed--quality trade-off that users can consider.

\subsection{DD-based cutting planes} \label{sec:exp-cuts}

We now illustrate how DD-suite can be used as a component of more sophisticated optimization algorithms and, in particular, how it enables reproducing results from the literature. Concretely, we recreate the 0--1 second-order conic (SOC) programming experiment of Castro, Cire, and Beck \cite{castro2022combinatorial}, who introduced the combinatorial cut-and-lift procedure that underlies our \emph{FlowCuts} and \emph{JointFlowCuts} generators and their strengthening; we additionally evaluate the target cuts of Tjandraatmadja and van Hoeve \cite{tjandraatmadja2019target} as a third cut family on the same benchmark. 

Concretely, we implement a DD-based cutting plane procedure and embed it into the Gurobi B\&B search through a user-cut callback at the root node (see Section \ref{sec:cutting-planes}). Following the setup of \cite{castro2022combinatorial}, for each SOC knapsack instance we build one relaxed DD per SOC constraint and, given a fractional point, iterate over the DDs to separate valid inequalities that cut off the point; enabling Gurobi's \texttt{PreCrush} parameter ensures that the resulting user cuts are applied to the presolved model. We evaluate three families of cuts available in DD-suite: combinatorial flow cuts (\emph{FlowCuts}), dual flow cuts derived from a max-flow formulation (\emph{JointFlowCuts}), and target cuts (\emph{TargetCuts}). For each family, we consider a plain variant and a strengthened variant (\emph{$+$Strengthening}) that lifts every generated inequality into a tighter one. Note that we omit the TargetCuts$+$Strengthening alternative since TargetCuts are known to yield facet-defining inequalities for this application \cite{tjandraatmadja2019target}. As a baseline, we run Gurobi with its default configuration.

Table \ref{tab:soc-cuts} reports aggregated results over the SOC-K benchmark for both the C++ and Python implementations, under a one-hour time limit. For each variant, we report the number of instances solved to optimality (\# Solv), the average final MIP gap over the unsolved instances (Gap), the average solution time over the solved instances (Time), the average number of B\&B nodes explored (\# Nodes), and the average number of cuts generated (\# Cuts). Averages are taken over the instances solved by both implementations to ensure a fair comparison. In each column and for each implementation, the best values are highlighted in bold.

\begin{table}[ht]
\centering
\caption{Aggregated results on the SOC-K benchmark.}
\label{tab:soc-cuts}
\setlength{\tabcolsep}{4pt}
\begin{tabular}{l rrrrr c rrrrr}
\toprule
 & \multicolumn{5}{c}{C++} & & \multicolumn{5}{c}{Python} \\
\cmidrule(lr){2-6} \cmidrule(lr){8-12}
Method & \# Solv & Gap & Time & \# Nodes & \# Cuts & & \# Solv & Gap & Time & \# Nodes & \# Cuts \\
\midrule
Gurobi
 & 81 & 0.68\% & 423.2 & 66{,}363 & --
 & & 81 & 0.67\% & 376.0 & 63{,}819 & -- \\
\midrule
FlowCuts
 & 83 & 0.73\% & 417.1 & 41{,}492 & 66
 & & 83 & 0.72\% & 393.5 & 40{,}847 & 66 \\
\quad $+$Strengthening
 & 84 & 0.62\% & 417.3 & 39{,}648 & 47
 & & 84 & 0.62\% & 418.3 & 39{,}658 & 46 \\
JointFlowCuts
 & 82 & 0.61\% & 333.1 & 30{,}883 & 174
 & & \textbf{86} & 0.61\% & 331.6 & 30{,}877 & 174 \\
\quad $+$Strengthening
 & 84 & \textbf{0.48\%} & 295.6 & 30{,}951 & 112
 & & 85 & \textbf{0.47\%} & 299.1 & 31{,}142 & 112 \\
TargetCuts
 & \textbf{85} & 0.54\% & \textbf{267.9} & \textbf{29{,}299} & 82
 & & \textbf{86} & 0.48\% & \textbf{294.6} & \textbf{28{,}127} & 83 \\
\bottomrule
\end{tabular}
\end{table}

The results show that all DD-based cut families improve upon the plain Gurobi baseline, which solves $81$ of the $90$ instances. Adding DD cuts increases the number of solved instances (up to $86$) and substantially reduces the size of the search tree: the number of B\&B nodes drops from roughly $66{,}000$ for Gurobi to about $30{,}000$ for the JointFlowCuts and TargetCuts variants, a reduction of more than $50\%$. The strengthening procedure further improves performance: for FlowCuts it reduces both the final gap (from $0.73\%$ to $0.62\%$) and the number of cuts needed (from $66$ to $47$), while for JointFlowCuts it attains the smallest final gaps of the whole table ($0.48\%$ in C++ and $0.47\%$ in Python). TargetCuts, in turn, solves the largest number of instances ($85$ in C++ and $86$ in Python) and yields the smallest average solution time ($268$ seconds in C++). Across the board, the JointFlowCuts and TargetCuts variants are the strongest, attaining the smallest gaps ($0.47\%$--$0.61\%$) and the smallest search trees. The C++ and Python implementations produce very similar results, as expected from their shared algorithmic core, with small differences attributable to solver-internal timing and tie-breaking.

We note that these results and conclusions are consistent with the ones obtained in \cite{castro2022combinatorial}, despite the latter using a specialized DD implementation in C++ with an older commercial solver (i.e., IBM ILOG CPLEX 12.9). Overall, this experiment demonstrates that DD-suite is not only a tool for building and querying DDs in isolation but also a practical foundation for developing competitive DD-based algorithms. 

\section{Conclusions and Future Work} \label{sec:conclusions}

This work introduced \textit{DD-suite}, a cross-platform, open-source code for building and manipulating decision diagrams (DDs) for discrete optimization. Unlike existing DD codes, which are typically tailored to a specific algorithm or application, DD-suite is designed to be clean, intuitive, and easily extended, providing the same modeling interface in both Python and C++. The suite supports the three standard DD types (exact, restricted, and relaxed), along with the DD reduction procedure, shortest-path routines for obtaining optimality bounds, a visualization tool, and extensive automated testing. Moreover, it includes thorough documentation, a support webpage, and ready-to-use examples for four combinatorial problems, as well as a SOC knapsack example used in our cutting plane experiments. Our goal is to lower the entry barrier to DD-based research, so that researchers can prototype and develop novel DD-based algorithms without reimplementing the underlying machinery.

Our numerical experiments support this design. The C++ and Python implementations build identical diagrams and bounds, with C++ being $5$--$6$ times faster and remaining within a small constant factor (about $2$ times) of \texttt{ddo} (i.e., a specialized, performance-oriented Rust framework). We further showed that DD-suite is easy to extend by adding an alternative grouping-based relaxation mechanism through a single overridden construction step, and that it serves as a practical building block for more sophisticated algorithms by implementing DD-based cutting planes that improve a state-of-the-art commercial solver on the SOC knapsack benchmark.

There are several promising directions for future work. A natural extension is to support additional DD types, such as zero-suppressed DDs, which are well suited to binary problems like the maximum independent set problem. Another is to incorporate additional DD construction mechanisms, such as iterative refinement or separation, which yield stronger relaxations and serve as building blocks for the column elimination algorithm. By releasing DD-suite as an open-source, well-documented, and extensively tested package, we hope to make DD-based optimization more accessible and to encourage its adoption by a broader research community.

\begin{acknowledgements}
 This research was partially supported by the National Center for Artificial Intelligence CENIA FB210017 (Basal ANID), and the \textit{Agencia Nacional de Investigación y Desarrollo} of Chile under grants ANID-FONDECYT-11230797 and ANID-FONDECYT-1261833.
\end{acknowledgements}

%
\section*{Conflict of interest}

The authors declare that they have no conflict of interest.

\bibliographystyle{spmpsci}      
\bibliography{references}   

\newpage
\begin{appendix}

\section{DD-suite Python code for our running example}\label{appendix:code}

\begin{lstlisting}[language=Python, 
    caption=Complete knapsack problem class implementation in Python,
    label=code:knapsackClass]
from SourceCode.Problems.AbstractProblemClass import AbstractProblem


class KnapsackProblem(AbstractProblem):

    def __init__(self, params: 'KnapsackStructure'):
        super().__init__(params.initial_state, params.variables)
        self.weights: list[int] = params.weights
        self.capacity: int = params.right_side_of_restrictions

    def transition_function(self, previous_state: 'State', variable_index: int, variable_value: int, scratch_state: list) -> bool:

        if variable_value == 0:
            scratch_state[0] = previous_state
            return True

        new_state: int = previous_state + self.weights[variable_index] * variable_value
        scratch_state[0] = new_state
        return new_state <= self.capacity

    def get_priority_for_discard_node(self, state: 'State') -> int:
        return -state

    def get_priority_for_merge_nodes(self, id: int, state: 'State') -> int:
        return -state

    def merge_operator(self, state_one: 'State', state_two: 'State') -> 'State':
        return min(state_one, state_two)

    def get_state_as_string(self, state: 'State') -> str:
        return str(state)

    def get_state_copy(self, state: 'State') -> 'State':
        return state
\end{lstlisting}

\begin{lstlisting}[language=Python, 
    caption=Complete DD construction and optimization code,
    label=code:knapsackMain]
from types import SimpleNamespace

from SourceCode.DD import DD
from Examples.KnapsackInstance.KnapsackProblem import KnapsackProblem
from SourceCode.GraphAlgorithms.ShortestLongestPath.ShortestLongestPath import ShortestLongestPath

## ============================
## SETUP PROBLEM AND PARAMETERS
## ============================

# Maximum width for the approximated DDs
width: int = 2

# Input data for the running example
params = SimpleNamespace(
    initial_state=0,
    variables=[('x_1', [0, 1]), ('x_2', [0, 1]), ('x_3', [0, 1]), ('x_4', [0, 1])],
    weights=[7, 5, 4, 1],
    right_side_of_restrictions=8,
)

# Objective coefficients and sense used later to optimize over the DDs
utility: list[int] = [4, 2, 5, 1]
objective: str = "max"

# Create the knapsack problem used to build the DD
knapsack_problem = KnapsackProblem(params)

## ============================
## CREATE EXACT DD
## ============================

dd_exact = DD(knapsack_problem)
dd_exact.create_decision_diagram()
dd_exact.reduce_decision_diagram()

# Print DD construction information
print("-- Exact DD information --")
print("\tConstruction time (sec): ", dd_exact.get_building_time())
print("\tReduction time (sec): ", dd_exact.get_reduction_time())
print("\tNumber of nodes: ", dd_exact.get_decision_diagram().get_node_count())
print("\tNumber of arcs: ", dd_exact.get_decision_diagram().get_arc_count())

## ============================
## CREATE RELAXED DD
## ============================

dd_relaxed = DD(knapsack_problem)
dd_relaxed.create_relax_priority_decision_diagram(max_width=width)
dd_relaxed.reduce_decision_diagram()

print("-- Relaxed DD information --")
print("\tConstruction time (sec): ", dd_relaxed.get_building_time())
print("\tReduction time (sec): ", dd_relaxed.get_reduction_time())
print("\tNumber of nodes: ", dd_relaxed.get_decision_diagram().get_node_count())
print("\tNumber of arcs: ", dd_relaxed.get_decision_diagram().get_arc_count())

## ============================
## CREATE RESTRICTED DD
## ============================

dd_restricted = DD(knapsack_problem)
dd_restricted.create_restricted_decision_diagram(max_width=width)
dd_restricted.reduce_decision_diagram()

print("-- Restricted DD information --")
print("\tConstruction time (sec): ", dd_restricted.get_building_time())
print("\tReduction time (sec): ", dd_restricted.get_reduction_time())
print("\tNumber of nodes: ", dd_restricted.get_decision_diagram().get_node_count())
print("\tNumber of arcs: ", dd_restricted.get_decision_diagram().get_arc_count())

# Export graph in .gml format
dd_exact.export_graph_file("knapsack_exact_dd_file")
dd_relaxed.export_graph_file("knapsack_relaxed_dd_file")
dd_restricted.export_graph_file("knapsack_restricted_dd_file")

## ============================
## OPTIMIZATION: PRIMAL AND DUAL BOUNDS
## ============================

# Optimal solution over the exact DD (longest path)
longest_path = ShortestLongestPath(dd_exact)
longest_path.set_parameters(utility, objective)
answer = longest_path.solve()

print("The optimal value is: ", answer.value)
print("Longest path: ", answer.path_print)
print("Solution time (sec): ", longest_path.get_time())

# Dual (upper) bound over the relaxed DD
longest_path_relaxed = ShortestLongestPath(dd_relaxed)
longest_path_relaxed.set_parameters(utility, objective)
answer_relaxed = longest_path_relaxed.solve()

print("The dual (upper) bound is: ", answer_relaxed.value)
print("Longest path over relaxed DD: ", answer_relaxed.path_print)
print("Solution time (sec): ", longest_path_relaxed.get_time())

# Primal (lower) bound over the restricted DD
longest_path_restricted = ShortestLongestPath(dd_restricted)
longest_path_restricted.set_parameters(utility, objective)
answer_restricted = longest_path_restricted.solve()

print("The primal (lower) bound is: ", answer_restricted.value)
print("Longest path over restricted DD: ", answer_restricted.path_print)
print("Solution time (sec): ", longest_path_restricted.get_time())
\end{lstlisting}

\section{Problems and data instances}\label{appendix:data-instances}

This appendix details, for each problem, the specific instance set used in the experiments of Section \ref{sec:experiments}. For each problem, we first state the classical integer programming (IP) formulation and then describe the instances themselves. In every case the recursive model assigns one variable per stage, so the DD horizon equals the number of variables (i.e., $\horizon = n$).

For knapsack, independent set, sequencing, and SOC knapsack we use established benchmarks from the literature; the sequencing instances are derived from a standard benchmark by the conversion described in Section \ref{appendix:data-sequencing}. For set cover, no standard benchmark fits our setting, so we generate the instances ourselves and report the generator below. Wherever random draws are involved, they use Python's \texttt{random} module (\texttt{numpy.random} for set cover) seeded with the value recorded in the file name, so every instance is reproducible from the listings. All the instances, together with the generation and conversion scripts, are available in the \texttt{DataInstances/} directory of the GitHub repository.

\subsection{Knapsack}\label{appendix:data-knapsack}

The \emph{$0$--$1$ knapsack problem} asks for a maximum-reward subset of $n$ items subject to a single capacity constraint. Letting $x_i = 1$ if item $i$ is placed in the knapsack, $c_i \in \Z$ its reward, $w_i \in \Z_{>0}$ its weight, and $b \in \Z_{>0}$ the capacity, the problem reads
\begin{align}\label{eq:ip-knapsack}
 \max\; & \sum_{i=1}^{n} c_i x_i \\
 \st\;  & \sum_{i=1}^{n} w_i x_i \leq b, \nonumber \\
        & \bx \in \B^n. \nonumber
\end{align}
The recursive formulation is that of Example \ref{example:recursiveKnapsack}: the state is the current load of the knapsack, and the resulting DD is binary.

We use the knapsack benchmark of \cite{pisinger2005hard}, which comprises the $30$ instances listed in Table \ref{tab:instances-knapsack} and combines two groups. The nine \emph{low-dimensional} instances (prefix \texttt{f}) are small on purpose: the number of items (i.e., binary variables) ranges from $4$ up to $23$, so that the exact DD can be constructed in full. The $21$ \emph{large-scale} instances (prefix \texttt{knapPI}) cover three instance types and seven sizes, from $100$ up to $10{,}000$ items, and are used to stress the restricted and relaxed constructions well beyond the low-dimensional regime.

\begin{table}[ht]
\centering
\caption{Knapsack instances \cite{pisinger2005hard}.}
\label{tab:instances-knapsack}
\footnotesize
\begin{tabular}{lll}
\toprule
\multicolumn{3}{c}{\emph{Low-dimensional}} \\
\midrule
\texttt{f1\_l-d\_kp\_10\_269} & \texttt{f2\_l-d\_kp\_20\_878} & \texttt{f3\_l-d\_kp\_4\_20} \\
\texttt{f4\_l-d\_kp\_4\_11} & \texttt{f6\_l-d\_kp\_10\_60} & \texttt{f7\_l-d\_kp\_7\_50} \\
\texttt{f8\_l-d\_kp\_23\_10000} & \texttt{f9\_l-d\_kp\_5\_80} & \texttt{f10\_l-d\_kp\_20\_879} \\
\midrule
\multicolumn{3}{c}{\emph{Large-scale}} \\
\midrule
\texttt{knapPI\_1\_100\_1000\_1} & \texttt{knapPI\_2\_100\_1000\_1} & \texttt{knapPI\_3\_100\_1000\_1} \\
\texttt{knapPI\_1\_200\_1000\_1} & \texttt{knapPI\_2\_200\_1000\_1} & \texttt{knapPI\_3\_200\_1000\_1} \\
\texttt{knapPI\_1\_500\_1000\_1} & \texttt{knapPI\_2\_500\_1000\_1} & \texttt{knapPI\_3\_500\_1000\_1} \\
\texttt{knapPI\_1\_1000\_1000\_1} & \texttt{knapPI\_2\_1000\_1000\_1} & \texttt{knapPI\_3\_1000\_1000\_1} \\
\texttt{knapPI\_1\_2000\_1000\_1} & \texttt{knapPI\_2\_2000\_1000\_1} & \texttt{knapPI\_3\_2000\_1000\_1} \\
\texttt{knapPI\_1\_5000\_1000\_1} & \texttt{knapPI\_2\_5000\_1000\_1} & \texttt{knapPI\_3\_5000\_1000\_1} \\
\texttt{knapPI\_1\_10000\_1000\_1} & \texttt{knapPI\_2\_10000\_1000\_1} & \texttt{knapPI\_3\_10000\_1000\_1} \\
\bottomrule
\end{tabular}
\end{table}

\subsection{Independent set}\label{appendix:data-independent-set}

Given an undirected graph $G = (V,E)$ with vertex rewards $c_i \in \Z$, the \emph{maximum weight independent set problem} looks for a set of pairwise non-adjacent vertices of maximum total reward. Using $x_i = 1$ to indicate that vertex $i$ is selected, the formulation is
\begin{align}\label{eq:ip-independent-set}
 \max\; & \sum_{i \in V} c_i x_i \\
 \st\;  & x_i + x_j \leq 1, \qquad \forall \{i,j\} \in E, \nonumber \\
        & \bx \in \B^{|V|}. \nonumber
\end{align}
The problem is also known as the \emph{maximum stable set} problem, and its unweighted version ($c_i = 1$ for all $i$) as the maximum independent set problem (MISP). It is equivalent to finding a maximum clique in the complement graph $\bar{G}$ and complementary to the minimum vertex cover problem, so care is needed when comparing objective values reported in the literature. 
In our recursive model, stage $t$ decides $x_t$ and the state records the set of vertices that are still eligible, which yields a binary DD with $n = |V|$ stages.

We use the DIMACS maximum-clique benchmark \cite{johnson1996dimacs}, a widely used set of graphs whose sizes range from a few dozen up to several hundred vertices. Since the maximum independent set of $G$ is a maximum clique of $\bar{G}$, the repository stores the \emph{complement} of each DIMACS graph and solves \eqref{eq:ip-independent-set} on it directly; for example, \texttt{MANN\_a9} has $45$ vertices and $918$ edges in DIMACS format and is stored with the $72$ edges of its complement. All instances are unweighted ($c_i = 1$), so the reported objective is the cardinality of the largest independent set, which coincides with the DIMACS maximum-clique number. We use the $64$ graphs listed in Table \ref{tab:instances-independent-set}, spanning the \texttt{MANN}, \texttt{brock}, \texttt{c-fat}, \texttt{hamming}, \texttt{johnson}, \texttt{keller}, \texttt{p\_hat}, \texttt{san}, and \texttt{sanr} families.

\begin{table}[ht]
\centering
\caption{Independent set instances (DIMACS \cite{johnson1996dimacs}).}
\label{tab:instances-independent-set}
\footnotesize
\begin{tabular}{llll}
\toprule
\texttt{brock200\_1} & \texttt{c-fat500-10} & \texttt{MANN\_a45} & \texttt{p\_hat700-3} \\
\texttt{brock200\_2} & \texttt{c-fat500-2} & \texttt{MANN\_a9} & \texttt{san1000} \\
\texttt{brock200\_3} & \texttt{c-fat500-5} & \texttt{p\_hat1000-1} & \texttt{san200\_0.7\_1} \\
\texttt{brock200\_4} & \texttt{hamming10-2} & \texttt{p\_hat1000-2} & \texttt{san200\_0.7\_2} \\
\texttt{brock400\_1} & \texttt{hamming10-4} & \texttt{p\_hat1000-3} & \texttt{san200\_0.9\_1} \\
\texttt{brock400\_2} & \texttt{hamming6-2} & \texttt{p\_hat1500-1} & \texttt{san200\_0.9\_2} \\
\texttt{brock400\_3} & \texttt{hamming6-4} & \texttt{p\_hat1500-2} & \texttt{san200\_0.9\_3} \\
\texttt{brock400\_4} & \texttt{hamming8-2} & \texttt{p\_hat1500-3} & \texttt{san400\_0.5\_1} \\
\texttt{brock800\_1} & \texttt{hamming8-4} & \texttt{p\_hat300-1} & \texttt{san400\_0.7\_1} \\
\texttt{brock800\_2} & \texttt{johnson16-2-4} & \texttt{p\_hat300-2} & \texttt{san400\_0.7\_2} \\
\texttt{brock800\_3} & \texttt{johnson32-2-4} & \texttt{p\_hat300-3} & \texttt{san400\_0.7\_3} \\
\texttt{brock800\_4} & \texttt{johnson8-2-4} & \texttt{p\_hat500-1} & \texttt{san400\_0.9\_1} \\
\texttt{c-fat200-1} & \texttt{johnson8-4-4} & \texttt{p\_hat500-2} & \texttt{sanr200\_0.7} \\
\texttt{c-fat200-2} & \texttt{keller4} & \texttt{p\_hat500-3} & \texttt{sanr200\_0.9} \\
\texttt{c-fat200-5} & \texttt{keller5} & \texttt{p\_hat700-1} & \texttt{sanr400\_0.5} \\
\texttt{c-fat500-1} & \texttt{MANN\_a27} & \texttt{p\_hat700-2} & \texttt{sanr400\_0.7} \\
\bottomrule
\end{tabular}
\end{table}

\subsection{Set cover}\label{appendix:data-set-cover}

Let $A \in \B^{m \times n}$ be a $0$--$1$ matrix whose $n$ columns index candidate sets and whose $m$ rows index the elements to be covered, and let $c_j \in \Z_{>0}$ be the cost of selecting column $j$. The \emph{set covering problem} selects a minimum-cost family of columns covering every row:
\begin{align}\label{eq:ip-set-cover}
 \min\; & \sum_{j=1}^{n} c_j x_j \\
 \st\;  & \sum_{j=1}^{n} A_{rj}\, x_j \geq 1, \qquad r = 1,\dots,m, \nonumber \\
        & \bx \in \B^n. \nonumber
\end{align}
When all costs are equal the problem reduces to minimum cardinality set cover, which is the case for every instance we use. Stage $t$ decides $x_t$ and the state records which rows remain uncovered, giving a binary DD with $n$ stages.

All the instances for this problem are generated by us, following the structured scheme of \cite{bergman2011manipulating} implemented in Listing \ref{code:gen-set-cover}. The key property of this scheme is that the matrix is \emph{banded}: row $i$ may only be covered by columns in the window $[i, \min(i+\beta, n))$, where $\beta$ is the bandwidth. Within that window the generator selects $k = \lfloor d\,n \rfloor$ columns uniformly without replacement (or the whole window when it contains fewer than $k$ columns) and sets those entries to one. The band structure is what makes these instances suitable for DDs: because a row can only interact with the $\beta$ columns that follow it, the number of rows that are simultaneously ``open'' at any stage is bounded, which keeps the exact DD width manageable and makes the effect of the maximum width $\maxwidth$ interpretable. All instances use unit costs, $c_j = 1$.

Instances follow the convention \texttt{set\_cover\_n\{n\}\_m\{m\}\_d\{d\}\_b\_w\{$\beta$\}\_seed\{S\}}. As summarized in Table \ref{tab:instances-set-cover}, they form two families, for a total of $384$ instances. The \emph{large} family fixes the bandwidth at $\beta = 165$ and sets $m = n - \beta + 1$, sweeping $n$ over the $16$ values $\{250, 500, \dots, 4000\}$ in steps of $250$. Here the density is not a free parameter: the generator sets $d = 75/n$, so that every row contains exactly $k = 75$ ones and the row density decreases as $n$ grows (this is why the file names show $d$ ranging from $0.30$ down to $0.02$). The \emph{small} family fixes $n = 100$ and $m = 78$ and instead sweeps the bandwidth over $\beta \in \{17, 20, \dots, 38\}$ and the density over $d \in \{0.05, 0.10, \dots, 0.50\}$, which isolates the effect of the bandwidth on the diagram size. Both families use seeds $S \in \{1,2,3,4\}$.

\begin{lstlisting}[language=Python, firstline=5,
    caption=Structured (banded) set cover instance generator.,
    label=code:gen-set-cover]
import numpy as np

def set_cover_matrix_generator(n: int, m: int, d: float, bw: int, seed: int=1) -> np.ndarray:
    """
    Generate an m x n matrix A for a structured Set Cover instance.

    Parameters:
    - n: number of columns (variables).
    - m: number of rows (constraints).
    - d: density of ones per row (a proportion between 0 and 1).
    - bw: bandwidth (how many columns ahead can be chosen).

    Returns:
    - A: binary matrix of size m x n.
    """
    A = np.zeros((m, n), dtype=int)
    k = int(d * n)  # number of ones per row
    np.random.seed(seed)

    for i in range(m):
        start = i
        end = min(i + bw, n)  # avoid exceeding n
        candidates = list(range(start, end))

        if len(candidates) == 0:
            continue  # there are no valid columns for this row

        selection = candidates if len(candidates) <= k else np.random.choice(candidates, size=k, replace=False)
        A[i, selection] = 1

    return A

n_options = [250, 500, 750, 1000, 1250, 1500, 1750, 2000, 2250, 2500, 2750, 3000, 3250, 3500, 3750, 4000]
d_options = [0.1]
b_w_options = [165]

for seed in range(1, 5):
    for n in n_options:
        for d in d_options:
            d = 75 / n  # fixes the number of ones per row at k = int(d * n) = 75
            for b_w in b_w_options:
                m = n - b_w + 1
                custom_folder = './Standard/'
                file_name = f'set_cover_n{n}_m{m}_d{round(d,2)}_b_w{b_w}_seed{seed}.txt'

                objective_weights = [1 for _ in range(int(n))]
                matrix_of_weight = set_cover_matrix_generator(n, m, d, b_w, seed)

                print(f"n: {n}, m: {m}, d: {d}, b_w: {b_w}")
                print(matrix_of_weight)

                with open(custom_folder + file_name, 'w') as file:
                    file.write(str(n) + '\n')
                    file.write(str(m) + '\n')
                    file.write(' '.join(map(str, objective_weights)) + '\n')
                    for i, row in enumerate(matrix_of_weight):
                        file.write(' '.join(map(str, row)))
                        if i != len(matrix_of_weight) - 1:
                            file.write('\n')
\end{lstlisting}

\begin{table}[ht]
\centering
\caption{Structured set cover instance families.}
\label{tab:instances-set-cover}
\footnotesize
\begin{tabular}{lccc}
\toprule
 & Large family & Small family \\
\midrule
$n$ (columns)      & $250, 500, \dots, 4000$ & $100$ \\
$m$ (rows)         & $n - \beta + 1$         & $78$ \\
Bandwidth $\beta$  & $165$                   & $17, 20, \dots, 38$ \\
Density $d$        & $75/n$                  & $0.05, 0.10, \dots, 0.50$ \\
Ones per row $k$   & $75$                    & $\lfloor d\,n \rfloor$ \\
Costs $c_j$        & $1$                     & $1$ \\
Seeds              & $1,\dots,4$             & $1,\dots,4$ \\
\midrule
\# configurations  & $16$                    & $80$ \\
\# instances       & $64$                    & $320$ \\
\bottomrule
\end{tabular}
\end{table}

\subsection{Sequencing}\label{appendix:data-sequencing}

We consider the single-machine sequencing problem with sequence-dependent setup times and total weighted completion time objective, denoted $1\,|\,s_{ij}\,|\,\sum_j w_j C_j$ in the standard scheduling notation. A set of $n$ jobs $J = \{1,\dots,n\}$ must be processed one at a time without preemption; job $j\in J$ has processing time $p_j$ and weight $w_j$, and a setup time $s_{ij}$ elapses when job $j$ is started immediately after job $i\in J$. A dummy job $0$ (the \emph{depot}) represents the initial state of the machine, so $s_{0j}$ is the setup required before the first job. The goal is to order the jobs so as to minimize $\sum_{j \in J} w_j C_j$, where $C_j$ is the completion time of job $j \in J$.

Introducing binary variables $y_{ij} = 1$ if job $j\in J$ is processed immediately after job $i\in J$ (with $i = 0$ for the first job) and continuous completion times $C_j \geq 0$, the formulation used in our experiments is
\begin{align}\label{eq:ip-sequencing}
 \min\; & \sum_{j \in J} w_j C_j \\
 \st\;  & \sum_{j \in J} y_{0j} = 1, \nonumber \\
        & y_{0j} + \sum_{i \in J,\, i \neq j} y_{ij} = 1, & \forall j \in J, \nonumber \\
        & \sum_{j \in J,\, j \neq i} y_{ij} \leq 1, & \forall i \in J, \nonumber \\
        & C_j \geq s_{0j} + p_j - M\,(1 - y_{0j}), & \forall j \in J, \nonumber \\
        & C_j \geq C_i + s_{ij} + p_j - M\,(1 - y_{ij}), & \forall i,j \in J,\, i \neq j, \nonumber \\
        & C_j \geq 0, \quad y_{ij} \in \B. \nonumber
\end{align}
The first three constraint families state that exactly one job starts the sequence, that every job has exactly one predecessor, and that every job has at most one successor. The big-$M$ constraints link the completion times to the selected sequence and simultaneously act as subtour elimination; we use $M = \sum_{j \in J} p_j + \sum_{j \in J} \max_{i \neq j} s_{ij}$.

In the recursive model, a state is the pair (set of already scheduled jobs, last scheduled job) and stage $t$ selects the job placed in the $t$-th position. Since the decision at each stage is a job index rather than a yes/no choice, the resulting diagram is a \emph{multivalued} DD (MDD) with $n$ stages and domain $\{1,\dots,n\}$, unlike the binary diagrams of the other four problems. This is also why in our implementation the cut separation of Section \ref{sec:exp-cuts} uses target cuts for this problem.

The setup time matrices are taken from the \emph{asymmetric traveling salesman problem} (ATSP) set of TSPLIB \cite{reinelt1991tsplib}, the standard library for this family of problems. 
%
Given an ATSP instance on $n$ cities with distance matrix $d$, we take city $0$ as the depot and cities $1,\dots,n-1$ as the jobs, so that a job sequence is an open Hamiltonian path rooted at the depot and the resulting instance has $n - 1$ jobs. The setup matrix is then $s_{0j} = d_{0,j+1}$ for the depot row and $s_{i+1,j} = d_{i+1,j+1}$ otherwise. The diagonal of a TSPLIB matrix carries a big-$M$ value (from $9999$ to $10^8$ depending on the instance) marking the arc a city cannot take to itself; under the mapping above it falls exactly on the entries $s_{i+1,i}$, which no feasible sequence ever traverses. 
%

TSPLIB supplies distances only, so processing times and weights are drawn following the scheme used in the sequencing literature cited by \cite{cire2013multivalued}, namely $p_j, w_j \sim \mathcal{U}\{1,\dots,10\}$, with the seed derived from the instance name so that the draw is reproducible. We keep the ATSP instances with fewer than $50$ cities, which yields the eight matrices reported in Table \ref{tab:sequencing-instances}, and generate $5$ replications of $(p, w)$ for each, giving $8 \times 5 = 40$ instances following the convention \texttt{sequencing\_atsp\_\{name\}\_n\{jobs\}\_seed\{S\}}. The setup time ranges reported in the table exclude the diagonal entries $s_{i+1,i}$, which no feasible sequence traverses and which carry the TSPLIB sentinel discussed above.

\begin{table}[htbp]
\centering
\caption{Sequencing instances derived from the TSPLIB ATSP set.}
\label{tab:sequencing-instances}
\begin{tabular}{lrrrr}
\toprule
ATSP instance & cities & jobs & $\min s_{ij}$ & $\max s_{ij}$ \\
\midrule
\texttt{br17}  & $17$ & $16$ & $0$  & $74$   \\
\texttt{ftv33} & $34$ & $33$ & $7$  & $332$  \\
\texttt{ftv35} & $36$ & $35$ & $7$  & $332$  \\
\texttt{ftv38} & $39$ & $38$ & $7$  & $332$  \\
\texttt{p43}   & $43$ & $42$ & $0$  & $5160$ \\
\texttt{ftv44} & $45$ & $44$ & $7$  & $332$  \\
\texttt{ftv47} & $48$ & $47$ & $7$  & $348$  \\
\texttt{ry48p} & $48$ & $47$ & $54$ & $2782$ \\
\bottomrule
\end{tabular}
\end{table}


\subsection{SOC knapsack}\label{appendix:data-soc}

The \emph{conic} or \emph{second-order cone (SOC) knapsack problem} replaces the linear capacity rows of \eqref{eq:ip-knapsack} by conic ones. Given rewards $c_i$, linear coefficients $a_{ji}$, conic coefficients $d_{ji}$, right-hand sides $b_j$, and a scalar $\Omega \geq 0$, the problem is
\begin{align}\label{eq:ip-soc-knapsack}
 \max\; & \sum_{i=1}^{n} c_i x_i \\
 \st\;  & \sum_{i=1}^{n} a_{ji}\, x_i \;+\; \Omega \sqrt{\sum_{i=1}^{n} d_{ji}^{2}\, x_i} \;\leq\; b_j, \qquad j = 1,\dots,m, \nonumber \\
        & \bx \in \B^n. \nonumber
\end{align}
Because $x_i \in \B$ implies $x_i = x_i^2$, the radicand equals $\|D_j \bx\|_2^2$ with $D_j = \mathrm{diag}(d_{j1},\dots,d_{jn})$, so each constraint can be written as $\bm{a}_j^\top \bx + \Omega \|D_j \bx\|_2 \leq b_j$ and is indeed a second-order cone constraint. 

The instances are the SOC knapsack benchmark from \cite{atamturk2010conic}, used exclusively in the cutting plane experiment of Section \ref{sec:exp-cuts}. They follow the convention \texttt{knapsack\_n\{n\}m\{m\}o\{$\Omega$\}b\{B\}\_\{S\}}, where $n \in \{100,125,150\}$ is the number of variables, $m \in \{10,20\}$ is the number of conic constraints, $\Omega \in \{1,3,5\}$ is the conic multiplier of \eqref{eq:ip-soc-knapsack}, $B = 5$ is fixed throughout the benchmark, and $S \in \{0,\dots,4\}$ indexes five replications. This gives $3 \times 2 \times 3 \times 5 = 90$ instances. 
%

The instance files store, in order: a line with $n$ and $m$; the value of $\Omega$; the reward vector $\bc$; the right-hand sides $b_1,\dots,b_m$; then $m$ lines with the linear coefficients $a_{ji}$; and finally $m$ lines with the conic coefficients $d_{ji}$.

\section{Verification across implementations}\label{appendix:exact-values}

To support the equivalence claim of Section \ref{sec:exp-efficiency}, we verify that: (i) the Python and C++ implementations build the same diagrams, (ii) the value of an exact DD is indeed the optimum of the problem, and (iii) the restricted and relaxed diagrams return valid primal and dual bounds.

The first check compares the two implementations across all possible configurations (i.e., combinations of instance, construction type, maximum width $\maxwidth$, reduction flag, and objective sense). For each such configuration, we compare three quantities: the optimization value returned by the diagram, the number of nodes of the diagram, and the number of arcs. Table \ref{tab:agreement-py-cpp} reports the results broken down by problem and construction type. Specifically, each cell represents the number of configurations on which the two implementations agree, out of the total number of configurations where both alternatives finish within the time limit.  The three quantities coincide in all $20{,}880$ configurations, which span $608$ instances, the five example problems, the four construction types, both the reduced and the non-reduced variants, and every maximum width used in Section \ref{sec:experiments}. 

\begin{table}[ht]
\centering
\caption{Python vs.\ C++ agreement, by problem and construction type.}
\label{tab:agreement-py-cpp}
\small
\setlength{\tabcolsep}{4pt}
\begin{tabular}{lccccc}
\toprule
Problem & Exact & Restricted & Relaxed (priority) & Relaxed (grouping) & Total \\
\midrule
Knapsack         & 45/45     & 319/319             & 318/318             & 295/295             & 977/977 \\
Independent set  & 36/36     & 766/766             & 754/754             & 725/725             & $2{,}281/2{,}281$ \\
Set cover        & 422/422   & $4{,}127/4{,}127$   & $4{,}125/4{,}125$   & $4{,}078/4{,}078$   & $12{,}752/12{,}752$ \\
SOC knapsack     & 180/180   & $1{,}080/1{,}080$   & $1{,}080/1{,}080$   & $1{,}080/1{,}080$   & $3{,}420/3{,}420$ \\
Sequencing       & 10/10     & 480/480             & 480/480             & 480/480             & $1{,}450/1{,}450$ \\
\midrule
Total            & 693/693   & $6{,}772/6{,}772$   & $6{,}757/6{,}757$   & $6{,}658/6{,}658$   & $20{,}880/20{,}880$ \\
\bottomrule
\end{tabular}
\end{table}

The second check validates the diagrams against an external reference: we compare the value of the (non-reduced) exact DD with the optimal value reported by Gurobi. This comparison is meaningful only for exact diagrams, since the restricted and relaxed constructions return bounds rather than optima. Table \ref{tab:exact-values} reports, for each problem, the number of instances on which the exact DD value matches Gurobi's optimum; it does so on every one of the $282$ instances. Two remarks are in order about the coverage of this comparison, which spans $282$ of the $518$ instances of these four problems. First, in the remaining instances, the exact DD exceeds the five-minute construction limit, which is precisely the regime that motivates the restricted and relaxed diagrams in Section \ref{sec:experiments}. Second, the converse also holds: in eight instances, the exact DD certifies an optimum that Gurobi fails to prove within its one-hour limit. These are one independent set instance and the five sequencing instances, where Gurobi terminates with optimality gaps between $57\%$ and $118\%$ while having already found the same solution, and two large knapsack instances, where Gurobi stops at its default relative optimality tolerance of $10^{-4}$ with an incumbent one unit below the true optimum. The SOC knapsack instances are excluded from this second comparison because each SOC knapsack DD models a single conic constraint and is used as a cut generator rather than as a standalone exact solver.

\begin{table}[tb]
\centering
\caption{Agreement of exact DD values with the optimum reported by Gurobi.}
\label{tab:exact-values}
\begin{tabular}{lcc}
\toprule
Problem & \# instances & DD $=$ Gurobi \\
\midrule
Knapsack         & 27  & 27 \\
Independent set  & 21  & 21 \\
Set cover     & 229 & 229 \\
Sequencing       & 5   & 5 \\
\midrule
Total            & 282 & 282 \\
\bottomrule
\end{tabular}
\end{table}

The third check verifies that the restricted and relaxed constructions return valid bounds. A restricted DD represents a subset of the feasible solutions, so its optimization bound is a primal bound that can never be better than the optimum; a relaxed DD represents a superset, so its optimization bound is a dual bound that can never be worse. For a maximization problem, we therefore require the restricted value to be at most the optimum and the relaxed value to be at least the optimum; the reverse holds for a minimization problem. We take as the reference the optimal value of the non-reduced exact DD when it is available; otherwise, we use the value reported by Gurobi for the instances where it closes the gap. This yields $490$ instances with a known optimum ($405$ certified by the exact DD and $85$ by Gurobi), making $17{,}589$ of the $21{,}297$ restricted and relaxed configurations checkable; the remainder correspond to instances whose optimum is unknown. The check spans the six maximum widths $\maxwidth$ of Section \ref{sec:experiments} and both the reduced and the non-reduced variants, and it is applied with a relative tolerance of $10^{-6}$, since the values are stored as floating-point numbers. Table \ref{tab:bounds-validity} reports, in each cell, the number of configurations whose bound holds over the number of configurations checked: every one of the $17{,}589$ bounds holds. As in the first check, only the C++ implementation is verified, since Table \ref{tab:agreement-py-cpp} already establishes that both implementations return the same values.

\begin{table}[ht]
\centering
\caption{Validity of the restricted and relaxed bounds.}
\label{tab:bounds-validity}
\small
\setlength{\tabcolsep}{4pt}
\begin{tabular}{lcccc}
\toprule
Problem & Restricted & Relaxed (priority) & Relaxed (grouping) & Total \\
\midrule
Knapsack         & 348/348           & 346/346           & 345/345           & $1{,}039/1{,}039$ \\
Independent set  & 492/492           & 484/484           & 490/490           & $1{,}466/1{,}466$ \\
Set cover        & $3{,}888/3{,}888$ & $3{,}888/3{,}888$ & $3{,}888/3{,}888$ & $11{,}664/11{,}664$ \\
SOC knapsack     & $1{,}080/1{,}080$ & $1{,}080/1{,}080$ & $1{,}080/1{,}080$ & $3{,}240/3{,}240$ \\
Sequencing       & 60/60             & 60/60             & 60/60             & 180/180 \\
\midrule
Total            & $5{,}868/5{,}868$ & $5{,}858/5{,}858$ & $5{,}863/5{,}863$ & $17{,}589/17{,}589$ \\
\bottomrule
\end{tabular}
\end{table}

\section{Gap of the relaxed bound: priority vs.\ grouping relaxation}\label{appendix:grouping-gaps}

\begin{figure}[ht]
    \centering
    \includegraphics[width=0.45\textwidth]{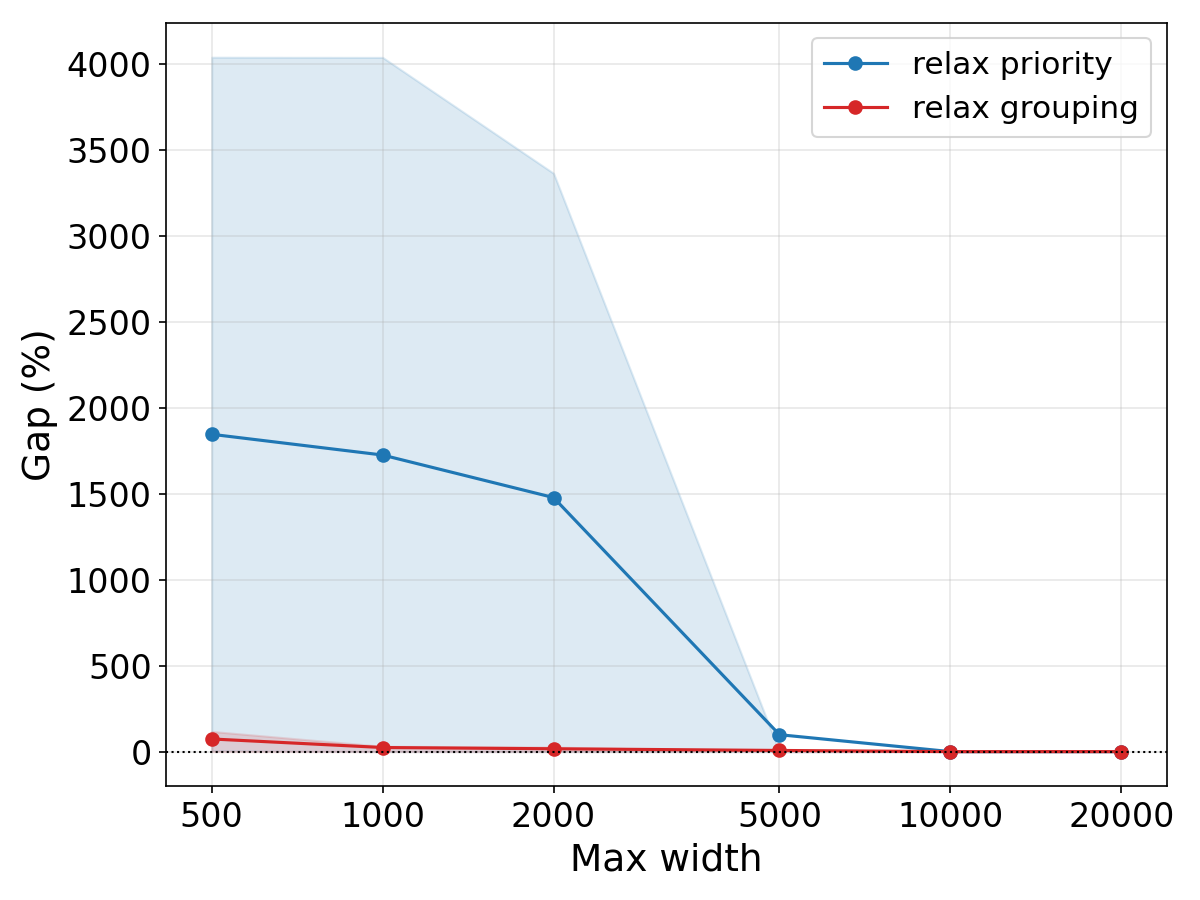}
    \includegraphics[width=0.45\textwidth]{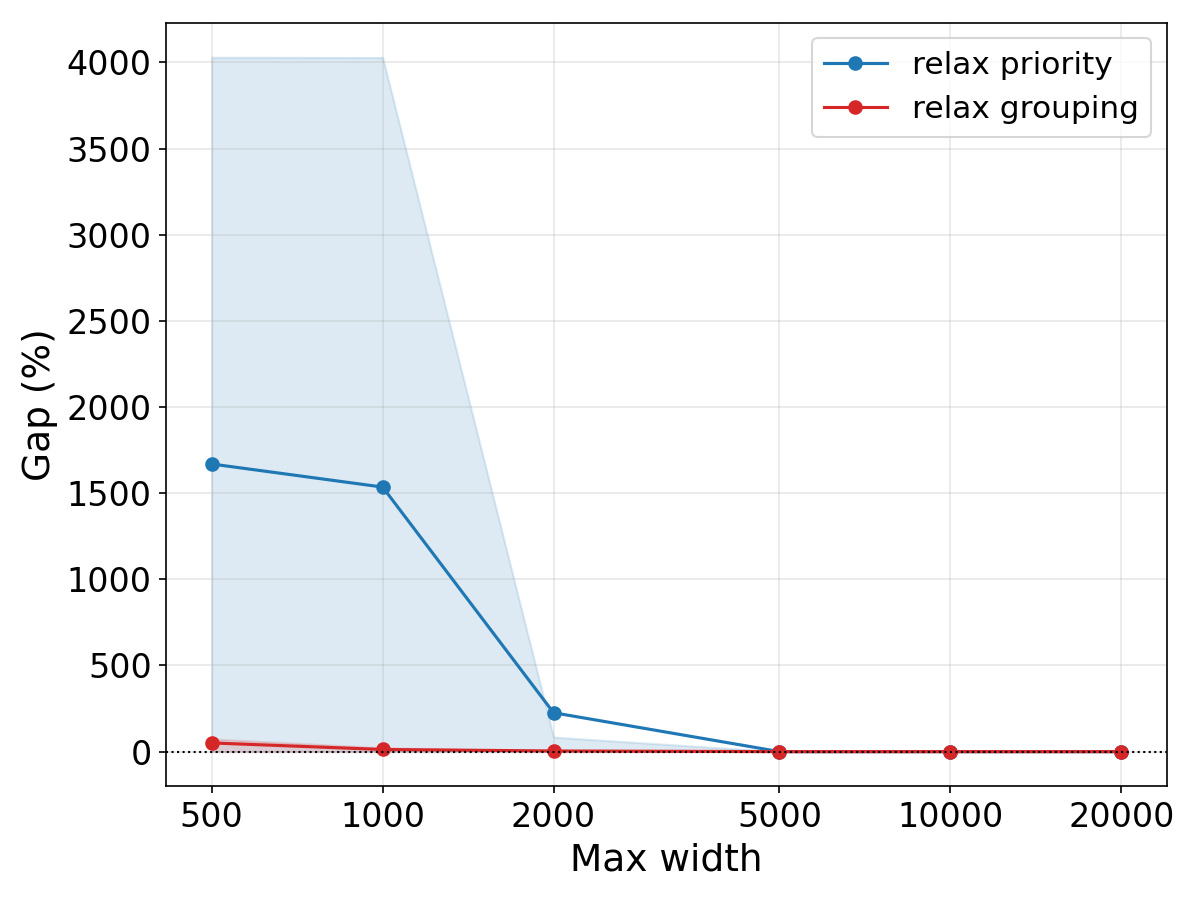}\\[0.5ex]
    \includegraphics[width=0.45\textwidth]{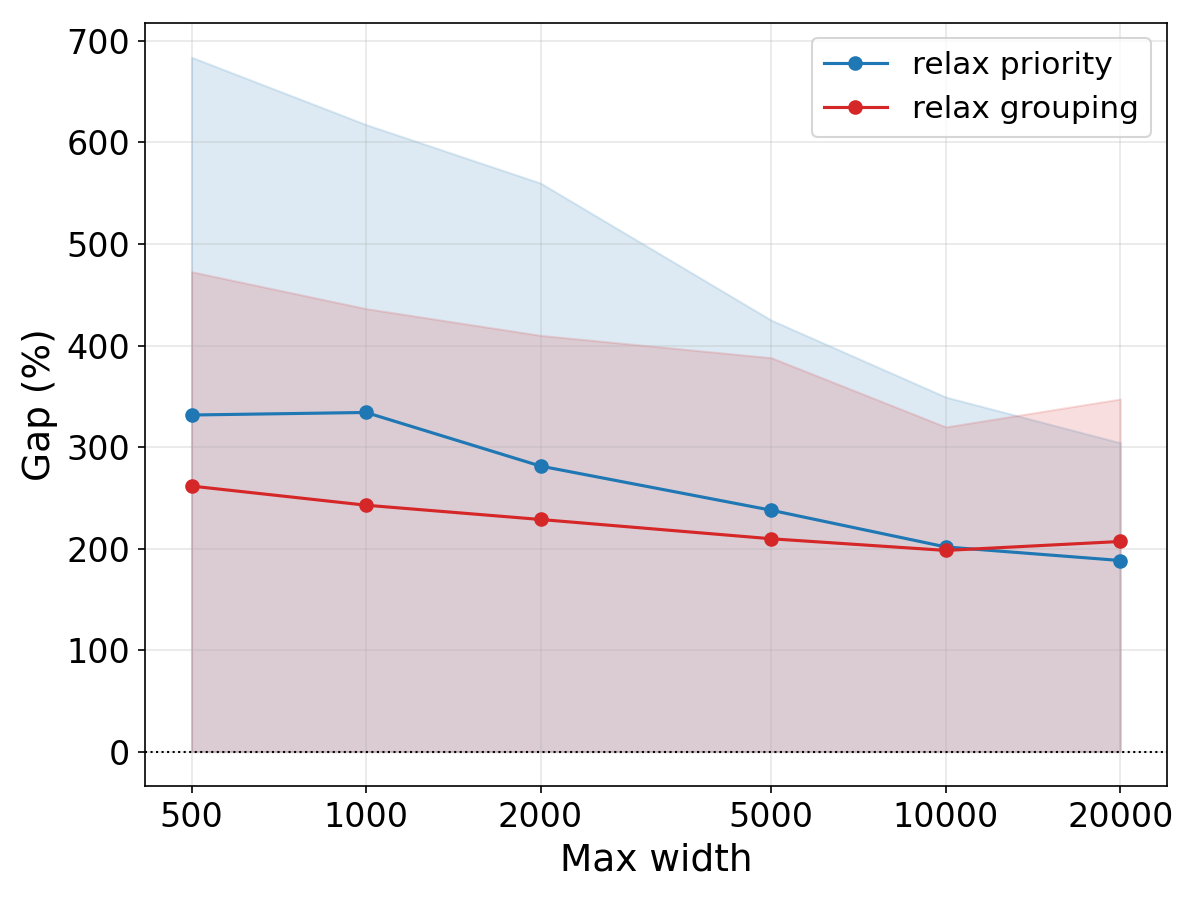}
    \includegraphics[width=0.45\textwidth]{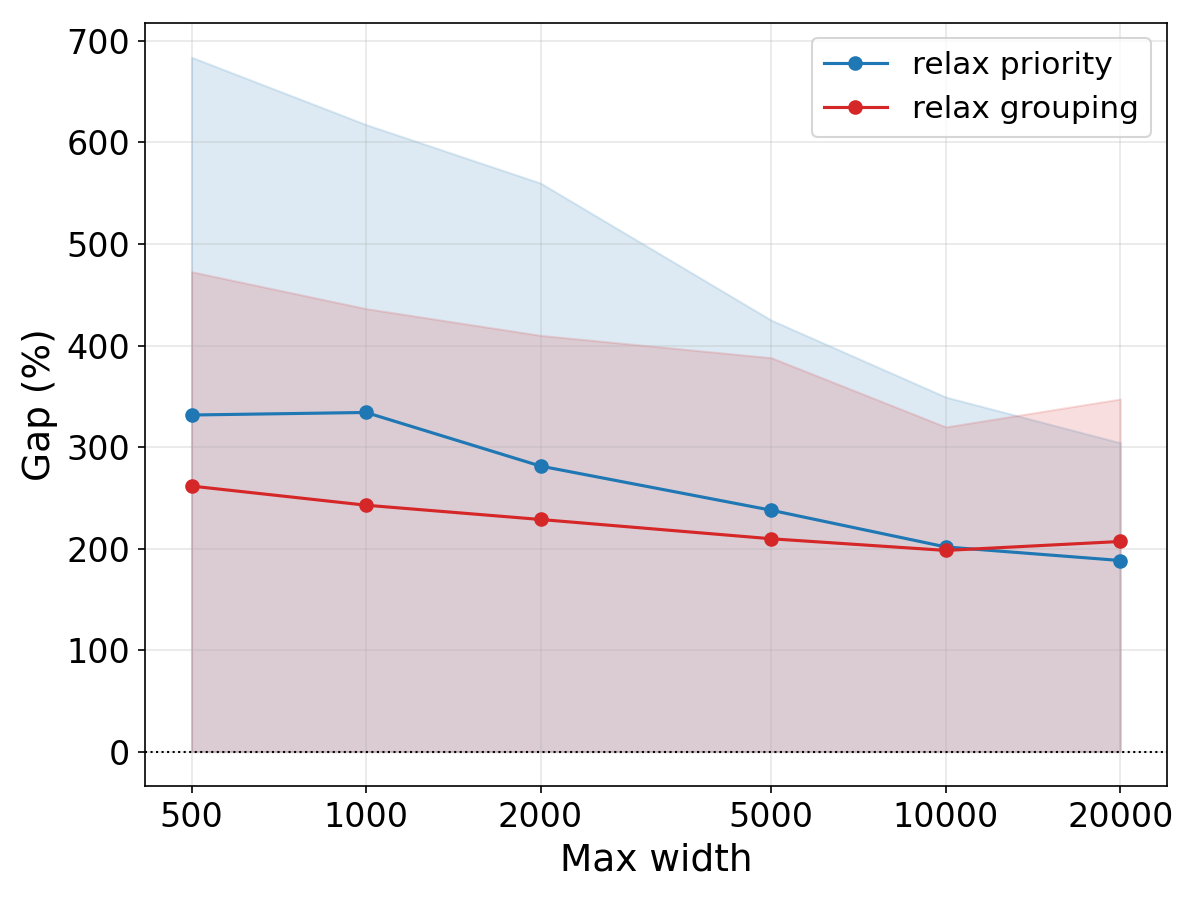}\\[0.5ex]
    \includegraphics[width=0.45\textwidth]{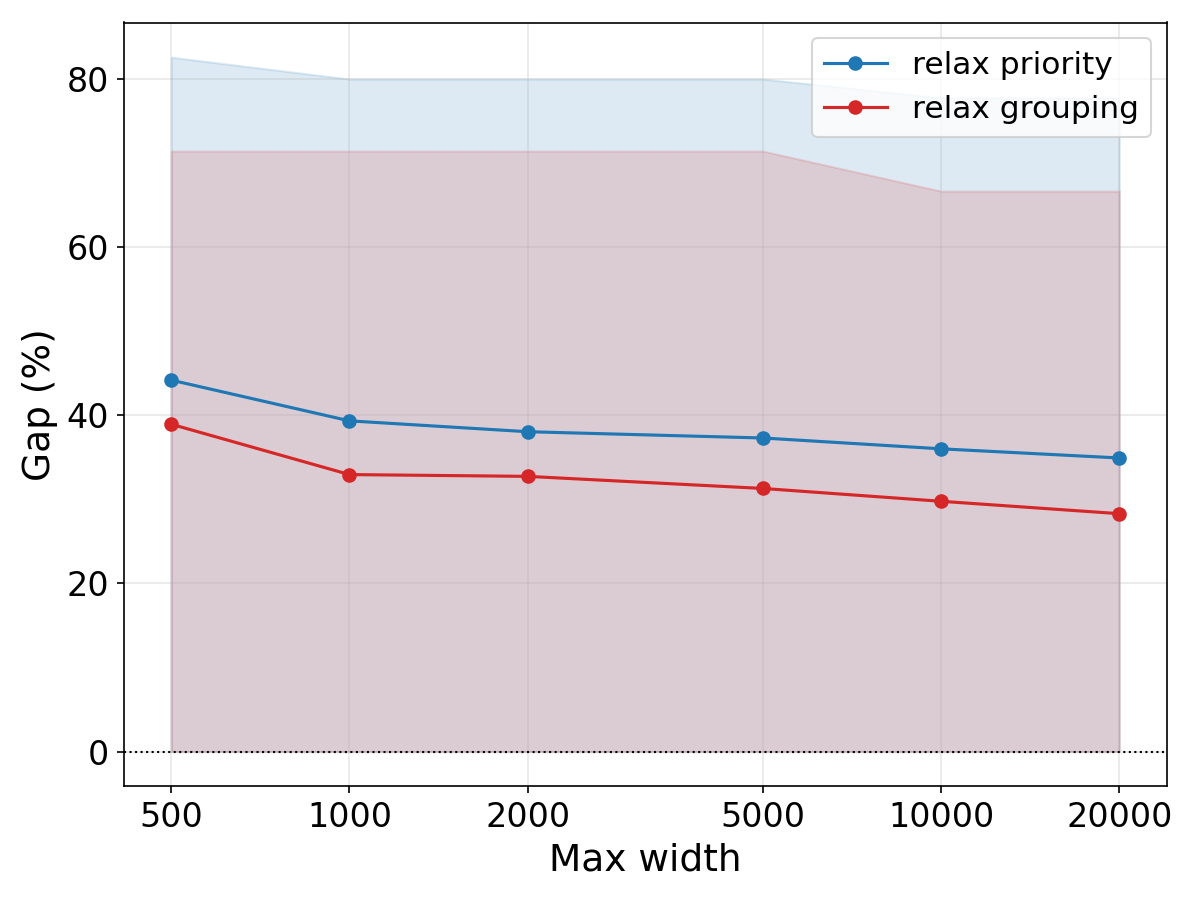}
    \includegraphics[width=0.45\textwidth]{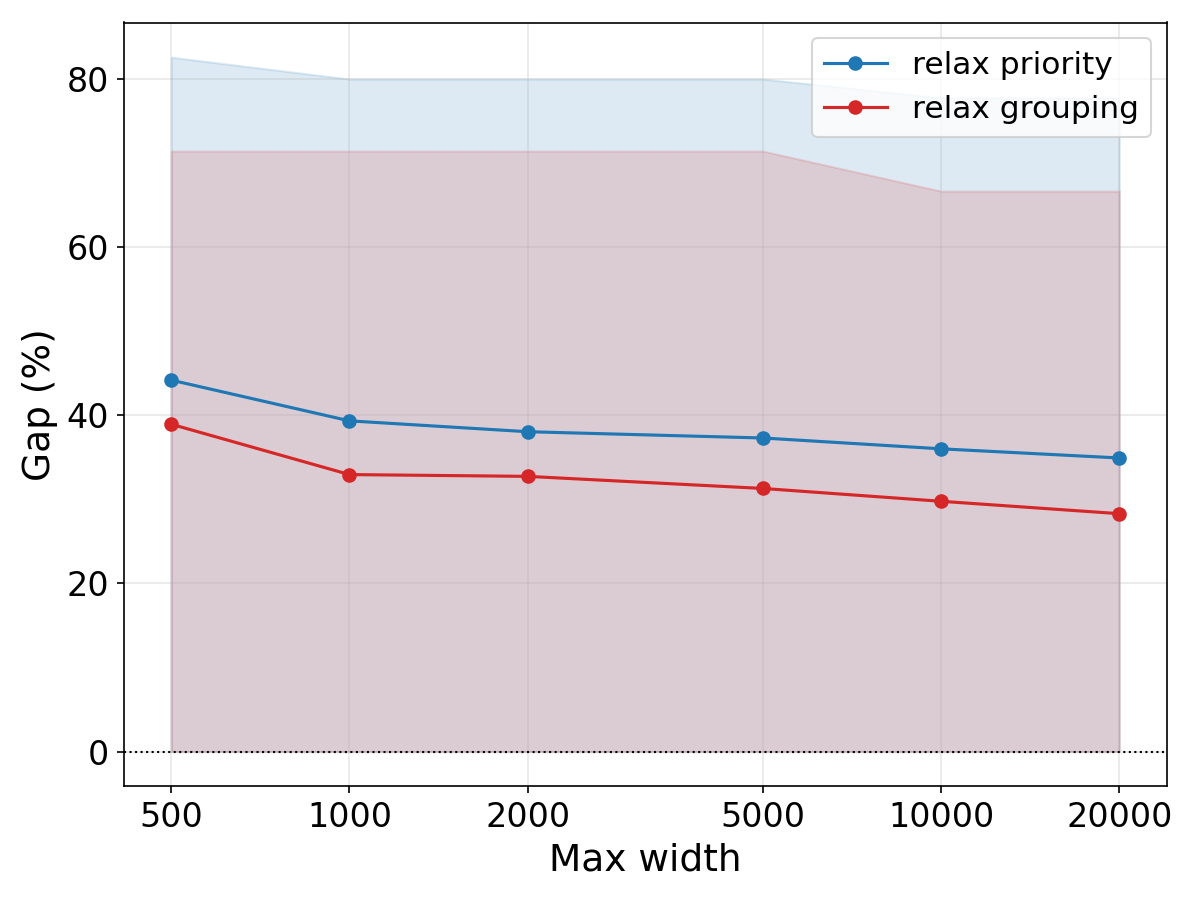}\\[0.5ex]
    \includegraphics[width=0.45\textwidth]{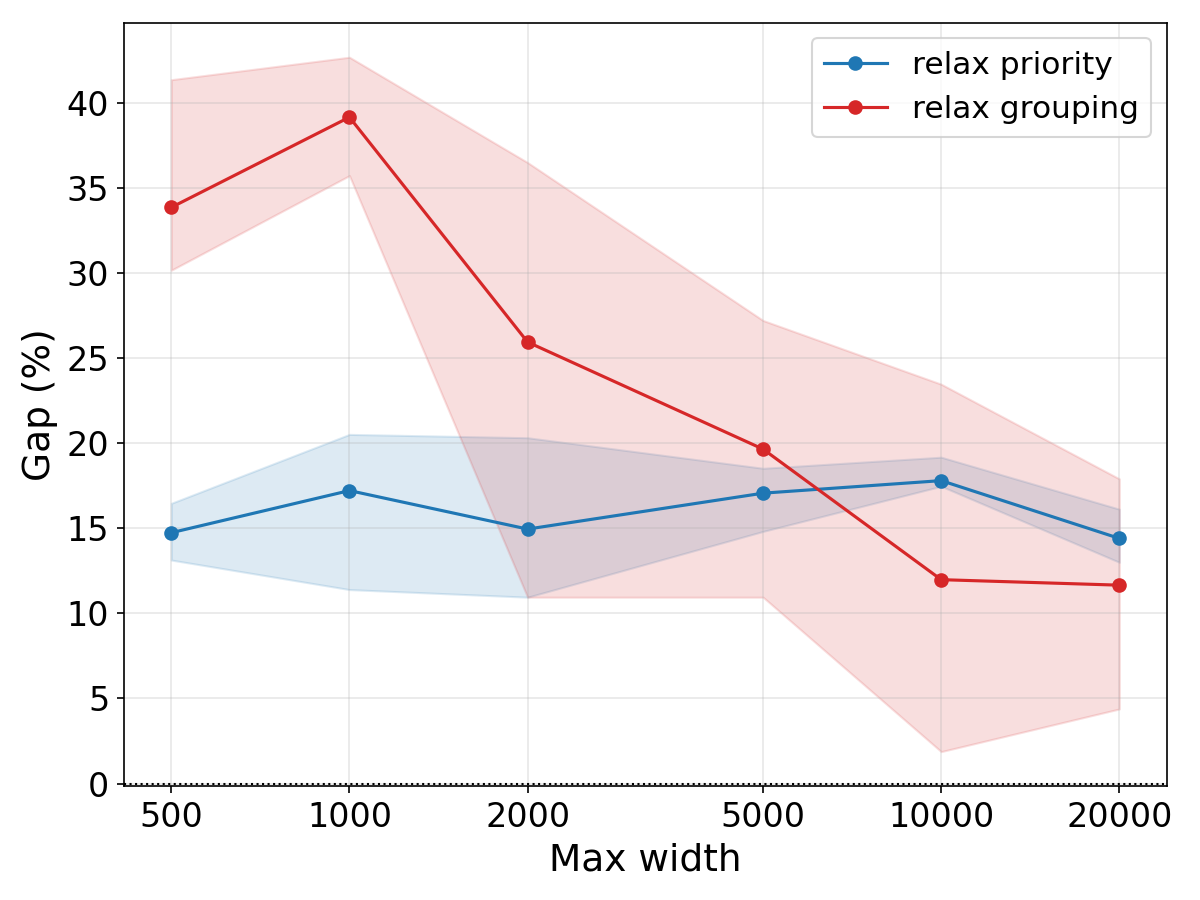}
    \includegraphics[width=0.45\textwidth]{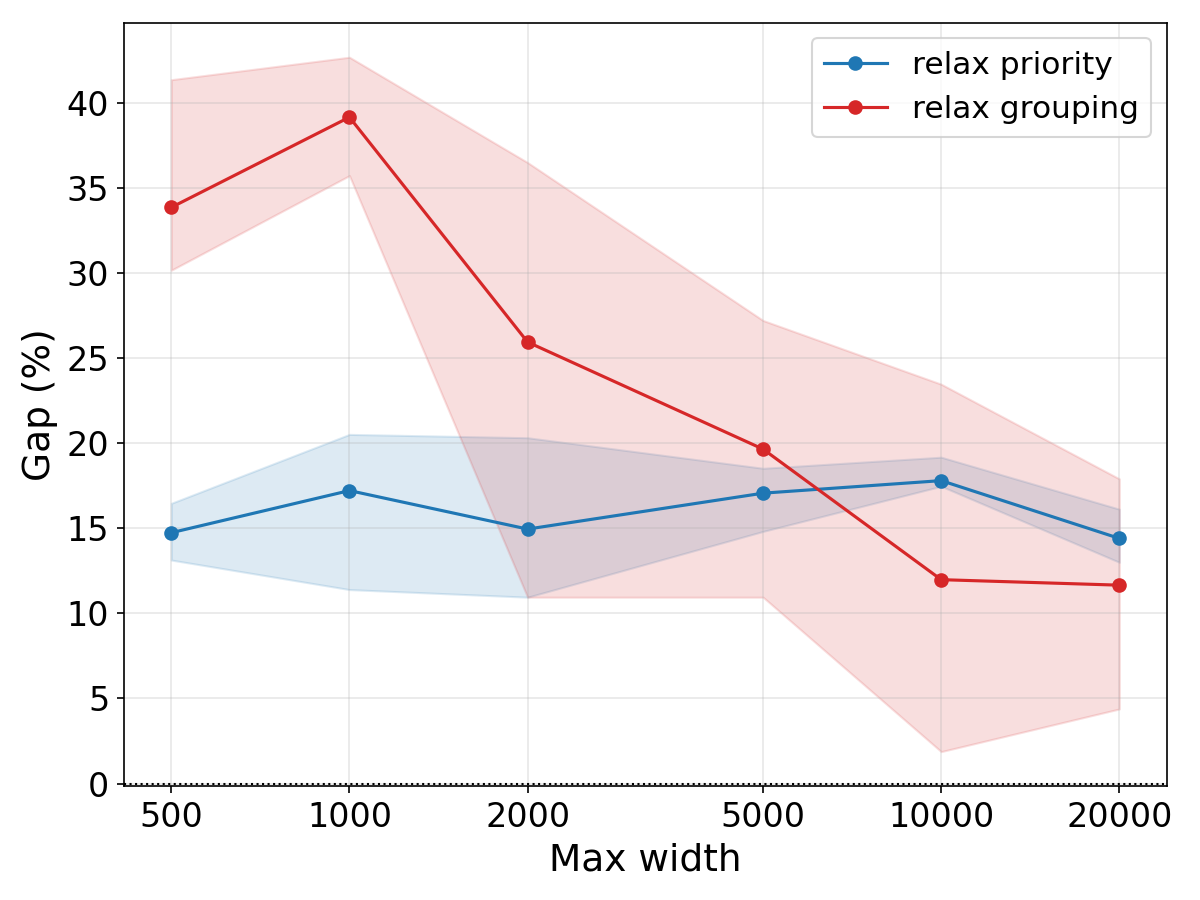}\\[0.5ex]
    \includegraphics[width=0.45\textwidth]{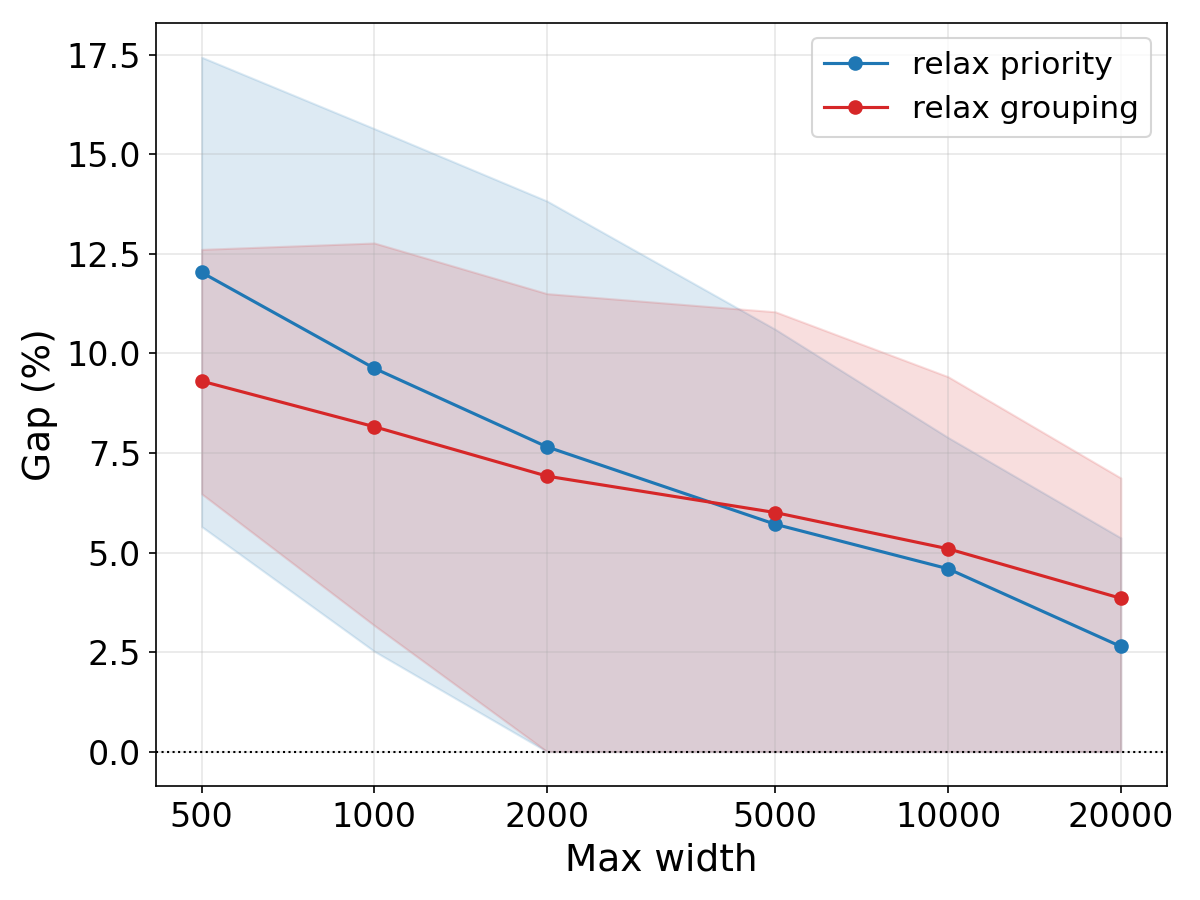}
    \includegraphics[width=0.45\textwidth]{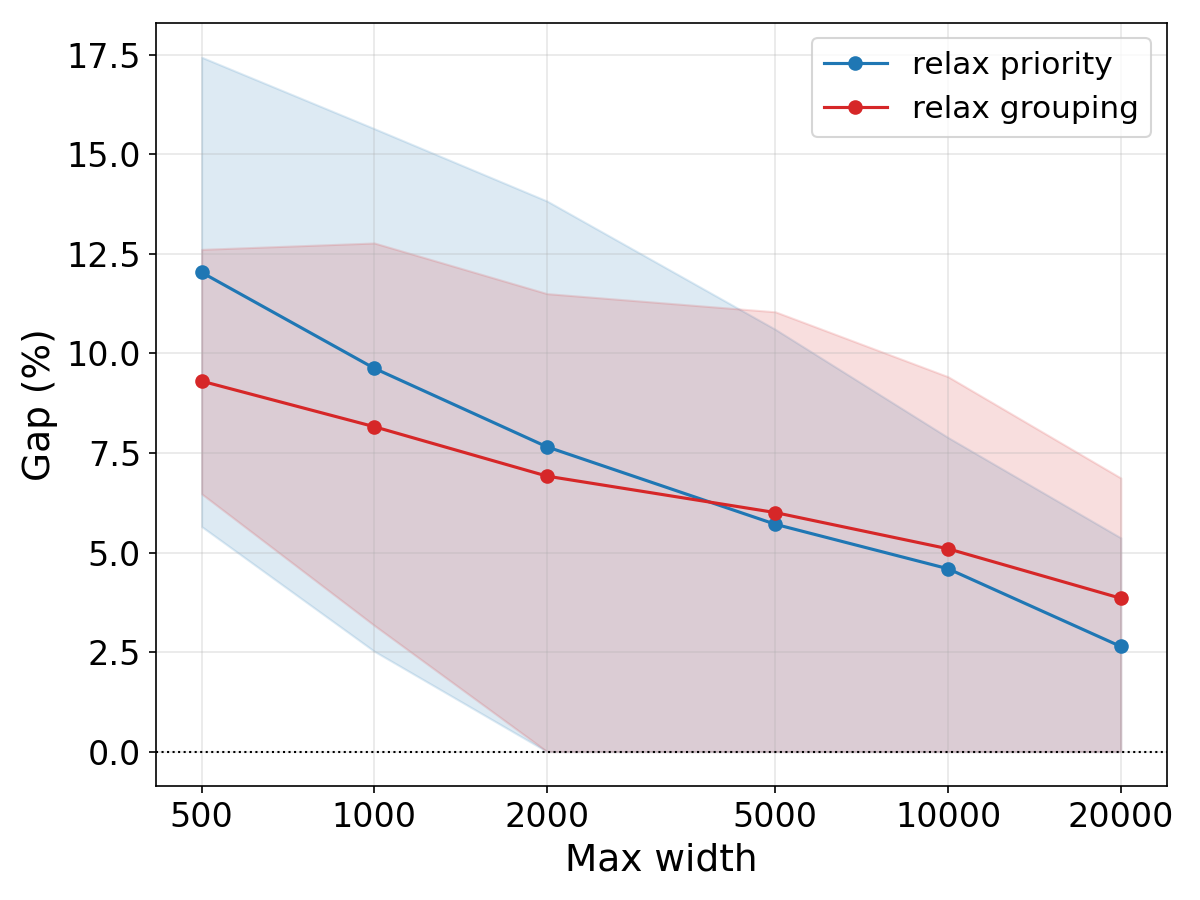}
    \caption{Relaxed-bound gap versus maximum width (all problems).}
    \label{fig:grouping-gaps-all}
\end{figure}

This appendix complements the extensibility experiment of Section \ref{sec:exp-grouping}. Figure \ref{fig:grouping-gaps-all} reports, for every example problem and for both the C++ and Python implementations, the optimality gap of the relaxed bound (relative to the exact optimum) as a function of the maximum width, for the two relaxation mechanisms (priority and grouping). Each row corresponds to a problem (from top to bottom: knapsack, independent set, set cover, sequencing, and SOC knapsack) and each column to an implementation (left: C++, right: Python). In each panel the line is the average gap after removing outliers with the IQR rule and the shaded band the interquartile range. These plots confirm the observation made in the main text: the two mechanisms yield comparable bound quality across all problems, with neither dominating the other.

\section{DD visualizations using yEd}\label{appendix:yEd_graphs}

\begin{figure}
    \centering
    \includegraphics[width=0.9\linewidth]{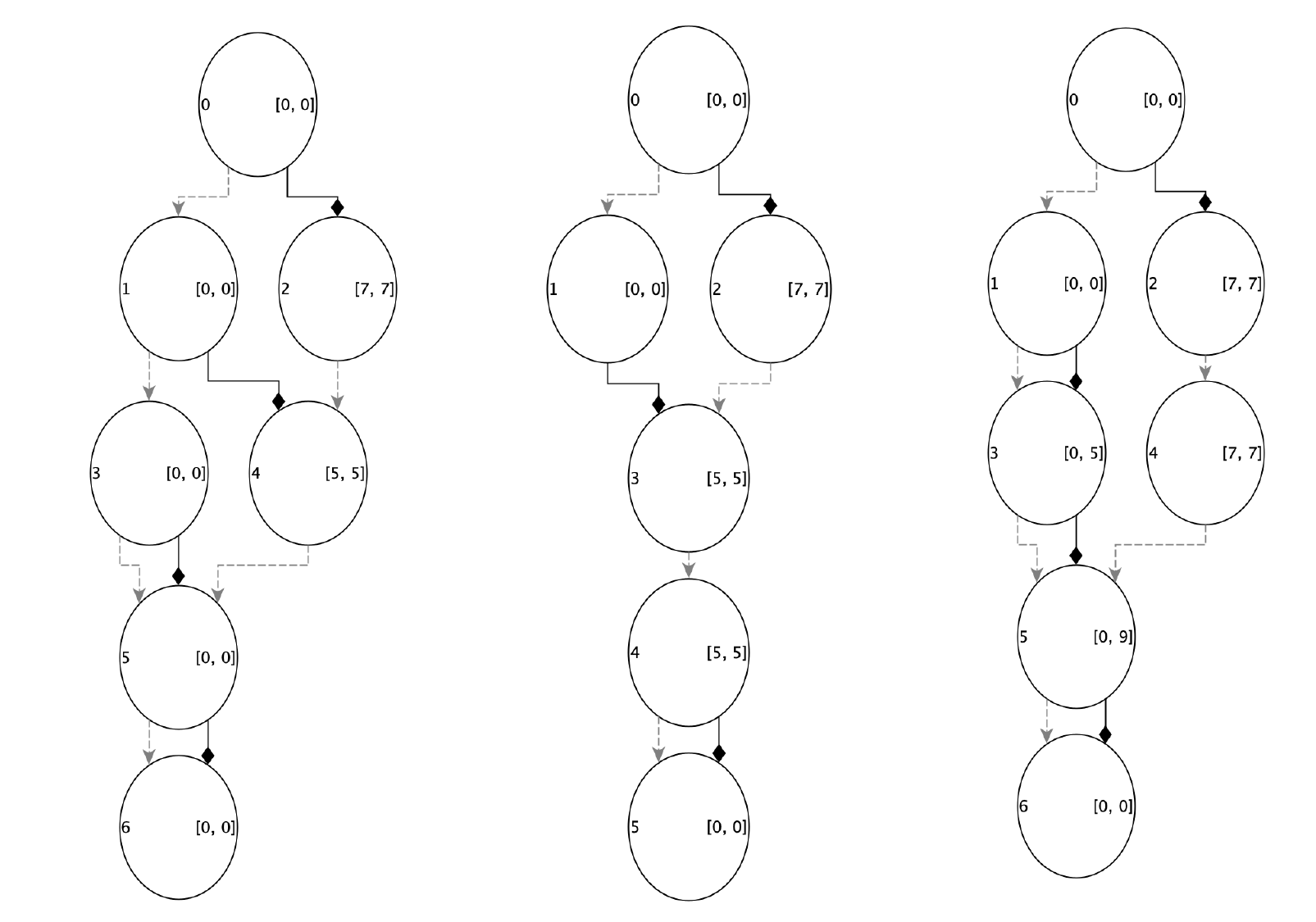}
    \caption{yEd visualizations of the three DD types for our running example.}
    \label{fig:yEd_DDgraphs}
\end{figure}

Figure \ref{fig:yEd_DDgraphs} shows yEd visualizations of the exact, restricted, and relaxed DDs for our running example. Note that the states shown for the exact and restricted DDs are not accurate, since these diagrams are reduced (i.e., nodes with different states are merged) and the merge operator is not used during the reduction procedure. Nonetheless, they encode the same set of feasible solutions as the diagrams in Figures \ref{fig:exactDDKnapsack} and \ref{fig:approxDDKnapsack}.

\end{appendix}

\end{document}